\documentclass[11pt]{amsart}
\usepackage{microtype}
\usepackage{hyperref}
\usepackage{color}
\usepackage{graphicx}
\usepackage[alphabetic]{amsrefs}

\newtheorem{thm}{Theorem}[section]
\newtheorem{lem}[thm]{Lemma}
\newtheorem{prop}[thm]{Proposition}

\newtheorem{conj}[thm]{Conjecture}

\newcommand{\newtheorembox}[2]{\newtheorem{x#1}[thm]{#2}
    \newenvironment{#1}{\begin{x#1}}{\qed\end{x#1}}}

\theoremstyle{definition}
\newtheorem{defn}[thm]{Definition}
\newtheorembox{ex}{Example}
\newtheorembox{rmk}{Remark}
\newtheorembox{const}{Construction}

\newcommand{\arxiv}[1]{\href{http://arxiv.org/abs/#1}{arXiv:#1}}

\newcommand{\CC}{\mathbb C}

\newcommand{\defi}[1]{\emph{#1}}
\newcommand{\degree}{\operatorname{deg}}
\renewcommand{\div}{\operatorname{div}}
\newcommand{\Gm}{\mathbb G_m}

\newcommand{\isom}{\cong}
\newcommand{\m}{\mathfrak m}
\newcommand{\mult}{\operatorname{mult}}

\newcommand{\PP}{\mathbb P}

\newcommand{\RR}{\mathbb R}
\newcommand{\Spec}{\operatorname{Spec}}

\newcommand{\val}{\operatorname{val}}
\newcommand{\X}{\mathfrak X}
\newcommand{\ZZ}{\mathbb Z}

\title{A specialization inequality for tropical complexes}
\author{Dustin Cartwright}
\address{Department of Mathematics \\ University of Tennessee \\ 227 Ayres Hall
         Knoxville, TN 37996}
\email{cartwright@utk.edu}

\begin{document}

\begin{abstract}
We prove a specialization inequality relating the dimension of the complete
linear series on a variety to the tropical complex of a regular semistable
degeneration. Our result extends Baker's specialization inequality to arbitrary
dimension.
\end{abstract}

\maketitle

\section{Introduction}

The specialization inequality for curves~\cite{baker} gives a bound for the
dimension of a linear system on a curve in terms of an analogous combinatorial
invariant on the dual graph of a degeneration. For higher-dimensional varieties,
understanding linear equivalence on the dual complex of a semistable
degeneration requires additional information beyond the dual complex, which can
be encoded in a tropical complex, as introduced in~\cite{cartwright-complexes}.
In this paper, we generalize the specialization inequality to varieties of
arbitrary dimension using tropical complexes.

Similar to the case of curves, our specialization inequality applies to a
regular, strictly semistable degeneration~$\X$ over a discrete valuation ring,
meaning that the special fiber of~$\X$ is a reduced union of smooth varieties,
with simple normal crossings. From $\X$, we can construct the dual
complex~$\Delta$, which is a regular $\Delta$-complex recording how the
components of the special fiber intersect. In addition, the tropical
complex~$\Delta$ records certain intersection numbers from $\X$, the details of
which will be recalled in Section~\ref{s:tropical-complexes}. In addition,
\cite{cartwright-complexes} introduced both a specialization map~$\rho$ from
divisors on the general fiber of~$\X$ to the tropical complex~$\Delta$ as well
as a compatible notion of linear equivalence for divisors on~$\Delta$. We define
$h^0$ of a divisor on~$\Delta$ as the fewest number of rational points such that
no linearly effective divisor contains all of these points (see Def.~\ref{d:h0}
for details), and then prove a specialization inequality:

\begin{thm}\label{t:specialization}
Let $\X$ be a regular strictly semistable degeneration of relative dimension~$n$
over a discrete valuation ring. Suppose that the locally closed
strata of dimension at most $n-2$ in $\X$ are affine and that $\X$ is robust in
dimensions $n-1$ and $n$. If $D$ is any divisor on the general fiber $X$
of~$\X$, and $\Delta$ is the tropical complex of~$\X$, then we have the
inequality:
\begin{equation*}
\dim H^0(X, \mathcal O(D)) \leq h^0(\Delta, \rho(D)).
\end{equation*}
\end{thm}

The second sentence of Theorem~\ref{t:specialization} places additional
requirements on the degeneration beyond the strict semistability, whose
definitions we now explain. A \defi{closed stratum} of dimension $n-k$ in~$\X$
is a connected component of the intersection of any $k+1$ components of the
special fiber. Each closed stratum contains a corresponding locally closed
stratum, which is formed by removing all lower-dimensional closed
strata, and so Theorem~\ref{t:specialization} requires these differences
to be affine varieties, when their dimension is $n-2$ or less. For the closed
strata of dimension~$n-1$ and $n$, Theorem~\ref{t:specialization} puts a weaker
condition of \defi{robustness}, which means that the union of the lower
dimensional strata form a big divisor, in the sense of birational geometry. See
Definition~\ref{d:robust-affine-strata} for details on both of these
definitions.

Although Theorem~\ref{t:specialization} does not apply directly to all regular
semistable degenerations because of these hypotheses, it does apply after
first modifying by a sequence of blow-ups, at least for projective
degenerations:
\begin{prop}\label{p:blow-ups}
Suppose that $\X$ is regular, strictly semistable degeneration, and all the
components of the special fiber $\X_0$ are projective.
Then, there exists a series of blow-ups, with centers contained in the special
fiber, which results in a regular, strictly semistable degeneration $\X'$ whose
locally closed strata are all affine.
\end{prop}

We now outline the structure of the proof of Theorem~\ref{t:specialization},
which helps to explain the necessity of the robustness and affine hypotheses.
Following the proof of the specialization inequality for curves, our
specialization inequality essentially follows from our definition of~$h^0$ given
a specialization map which preserves linear equivalence, effectivity of
divisors, and point containment. Preservation of linear equivalence was proved
in~\cite{cartwright-complexes}. However, for an effective divisor~$D$, the
specialization,~$\rho(D)$, is not always effective, so we introduce a refined
specialization in Section~\ref{s:refined}, which is effective and is linearly
equivalent to~$\rho(D)$. Note that, unlike $\rho(D)$, the refined specialization
is not necessarily supported on the $(n-1)$-dimensional simplices of~$\Delta$,
and so our specialization theorem necessarily requires a framework that includes
divisors which intersect the interior of the $n$-dimensional simplices, as
explained in Section~\ref{s:tropical-complexes}.

Although our proof does not use this technology explicitly, the refined
specialization of a divisor~$D$, and its relationship to point containment, can
be understood in terms of the projection of $D$ to the skeleton of the
Berkovich analytification. More specifically, by results going back to
Berkovich~\cite{berkovich}, the dual complex~$\Delta$ embeds into, and is a
strong deformation retract of, the analytification of the general fiber of~$\X$.
Restricting to effective divisors in order to avoid cancellation, we have:
\begin{prop}\label{p:analytification} If $D$ is an effective divisor in the
general fiber of~$\X$, then the projection of the analytification of $D$ to the
skeleton defined by $\X$ is the union of the refined specialization of $D$
together with a finite set of polyhedra, all of dimensions at most $n-2$.
\end{prop}

While the projection of the analytification always preserves point containment,
Proposition~\ref{p:analytification} shows that the same can fail for the
refined specialization whenever the projection has components of codimension~2
or greater. The purpose of the robustness and affine hypotheses in
Theorem~\ref{t:specialization} is to guarantee that we can lift points from the
tropical complex to the algebraic variety in such a way that containing divisors
project to a set of codimension~1 on~$\Delta$. In the case of
curves, the set of codimension~2 in Proposition~\ref{p:analytification} is
necessarily empty, which is why the specialization inequality for curves
required no hypotheses beyond a strictly semistable degeneration.

As an application of Theorem~\ref{t:specialization}, we give an example of a
tropical complex which does not lift to any algebraic variety.

\begin{thm}\label{t:no-lift}
There exists a $2$-dimensional tropical complex~$\Delta$ which is not the
tropical complex of any regular semistable degeneration.
\end{thm}

\noindent The underlying $\Delta$-complex of the example in
Theorem~\ref{t:no-lift} is a triangulation of the product of a cycle and an
interval. Therefore, the dual complex is realizable from a degeneration of an
algebraic surface, such as the product of an elliptic curve with a projective
line. However, in Theorem~\ref{t:no-lift} the structure constants of the
tropical complex~$\Delta$ are chosen so that for any divisor $D$, $h^0(\Delta,
nD)$ grows at most linearly in~$n$, which combined with the existence of an
ample divisor on any smooth, proper surface, shows that $\Delta$ can't be the
tropical complex of a projective surface.

An alternative approach to a specialization inequality has been developed in
unpublished work of Eric Katz and June Huh. They work not with degenerations,
but with tropicalizations of surfaces embedded in $\mathbb G_m^N$ over a field
with trivial valuation. Moreover, their proof involves choosing linearly
equivalent divisors passing through the lowest dimensional toric strata, as
opposed to the proof of Theorem~\ref{t:specialization} in which points are
lifted from the $n$-dimensional locally closed strata. Therefore, the two
approaches should give distinct and possibly complementary bounds on the
dimension of linear series.

The rest of this paper is organized as follows.
Section~\ref{s:tropical-complexes} recalls the definition of tropical
complexes from~\cite{cartwright-complexes} as well as their main properties. In
Section~\ref{s:linear-series}, we define the invariant~$h^0$ from
Theorem~\ref{t:specialization} and look at some examples and applications.
Section~\ref{s:subdivisions} proves compatibility of tropical complexes under
ramified base changes followed by toroidal resolutions of singularities, which
are the tools used to give a refined specialization in
Section~\ref{s:refined}. Section~\ref{s:proof} has the proof of the
specialization inequality.

\subsection*{Acknowledgments}

Throughout this project, I've benefited from my conversations with Matt Baker,
Spencer Backman, Alex Fink, Christian Haase, Paul Hacking, June Huh, Eric Katz,
Madhusudan Manjunath, Farbod Shokrieh, Bernd Sturmfels, Yu-jong Tzeng, and
Josephine Yu. I'd especially like to thank Sam Payne for his many insightful
suggestions and thoughtful comments on an earlier draft of this paper. I was
supported by the National Science Foundation award number DMS-1103856 and
National Security Agency Young Investigator Grant H98230-16-1-0019.

\section{Tropical complexes}\label{s:tropical-complexes}

Let $\X$ be a regular, strictly semistable degeneration over a discrete
valuation ring $R$. Specifically, we mean that $\X$ is a regular scheme, flat
and proper over $\Spec R$, and such that the fiber over the closed point of
$\Spec R$
is a reduced simple normal crossing divisor. The closed fiber is called the
\defi{special fiber} and denoted $\X_0$. We
also assume that the residue field of $R$ is algebraically closed. Under these
assumptions, we will call $\X$ a \defi{degeneration} and we let $n$ denote the
dimension of the general fiber of $\X$, which is also the dimension of each
component of the special fiber. The \defi{general fiber} of $\X$, denoted $X$, is the
fiber over the generic point of $\Spec R$.

We first recall the construction of the dual complex of $\X$. The dual complex
is a regular $\Delta$-complex (in the sense of \cite{hatcher}*{Sec. 2.1}) that
consists of a vertex for each irreducible component of the special fiber $\X_0$.
In addition, for each set of $k+1$ of the irreducible components of $\X_0$,
their intersection is a disjoint union of smooth $(n-k)$-dimensional varieties,
and we have a $k$-dimensional simplex for each of these components. We write
$C_s$ for the smooth variety corresponding to a $(k+1)$-dimensional variety, and
we have that the simplices are attached such that $C_s$ is a subvariety of
$C_{s'}$ if and only if $s'$ is contained in $s$. See the beginning of Section 2
of \cite{cartwright-complexes} for full details.

The \defi{weak tropical complex} of $\X$ is the pair of the dual complex and a
function $\alpha$ from the pairs of a vertex $v$ and a ridge $r$ which contains
$v$. Here, by \defi{ridge}, we mean an $(n-1)$-dimensional simplex of the dual
complex of $\X$. Thus, the corresponding component $C_r$ is a curve which is
contained in $C_v$, which is a component of $\X_0$ and a divisor on $\X$. We
define $\alpha(v,r) = - \degree C_v \cdot C_r$, where $\deg C_v \cdot C_r$
denotes the degree of the intersection product on $\X$. Note that the
intersection is necessarily non-transverse, because $v$ is a vertex of $r$. In
the transverse case, the intersection number would be redundant because it
number of points of $C_v \cap C_r$, which are in bijection with the
$n$-dimensional simplices of the dual complex containing both $v$ and $r$. The
weak tropical complex of a degeneration satisfies the following definition by
\cite{cartwright-complexes}*{Prop.~2.6}.

\begin{defn}[Def.~2.5 in \cite{cartwright-complexes}]\label{d:weak-tropical-complex}
An \defi{$n$-dimensional weak tropical complex} $\Delta$ is a pair $(S,\alpha)$,
where $S$ a finite, connected, regular $\Delta$-complex,
whose simplices all have dimension at most $n$, and
$\alpha$ is a function from pairs of a vertex $v$ of $S$ and a ridge~$r$
containing $v$ such that for every ridge $r$,
\begin{equation}\label{eq:ridge-identity}
\sum_{v \in r_0} \alpha(v,r) = \degree r,
\end{equation}
where $r_0$ denotes the vertices of $r$ and $\degree r$ is the number
of $n$-dimensional simplices containing $r$. The values $\alpha(v,r)$ are
sometimes referred to as the \defi{structure constants} of the weak tropical
complex. We refer to the $0$-, $1$-,
$(n-1)$-, and $n$-dimensional simplices of $S$ as the \defi{vertices},
\defi{edges}, \defi{ridges}, and \defi{facets}, respectively,
of~$\Delta$.
\end{defn}

Definition~\ref{d:weak-tropical-complex} is more specialized than the one
in~\cite{cartwright-complexes} in that we assume that the underlying
$\Delta$-complex is regular, meaning that distinct faces of a single simplex are
not identified with each other. The dual complex of a simple normal crossing
divisor, and thus of a degeneration, is always a regular $\Delta$-complex, and
thus regular $\Delta$-complexes are sufficient for the weak tropical complexes
which appear in this paper. Assuming that the $\Delta$-complex is regular
simplifies the exposition at several points. 

We now recall the definitions of locally closed strata and robustness, which
appear in the statement of Theorem~\ref{t:specialization}. Recall that a divisor
$D$ on a smooth variety is called \defi{big} if the complete linear series on
some multiple of $D$ defines a birational map onto its image in projective
space~\cite{lazarsfeld}*{Sec. 2.2}.

\begin{defn}\label{d:robust-affine-strata}
Let $\X$ be a degeneration and let $\Delta$ be its weak tropical complex. For
any $(n-k)$-dimensional simplex $s$ of $\Delta$, let $D_s$ denote the union
$\bigcup_{s'} C_{s'}$, where $s'$ ranges over $(n-k+1)$-dimensional simplices
containing~$s'$. We define $C_s \setminus D_s$ to be the \defi{locally closed
stratum} of dimension $k$, corresponding to $s$. We also consider $D_s$ as a
reduced divisor on $C_s$, and we say that $\X$ is \defi{robust in dimension $k$}
if, for any $(n-k)$-dimensional simplex~$s$, $D_s$ is a big divisor on $C_s$.
\end{defn}

\begin{prop}\label{p:robust-dim-1}
Let $\X$ be a degeneration, and $\Delta$ its dual complex. Then, the following
are equivalent:
\begin{enumerate}
\item $\X$ is robust in dimension 1.
\item The locally closed strata of dimension 1 in $\X$ are affine.
\item Every $(n-1)$-dimensional simplex of $\Delta$ is contained in some
$n$-dimensional simplex.
\end{enumerate}
\end{prop}

\begin{proof}
Let $r$ be an $(n-1)$-dimensional simplex of $\Delta$, and thus a ridge. Then,
$C_r$ is a curve, and $D_r$ is the union of points corresponding to facets
containing~$r$, and thus is non-trivial if and only if there are any such
facets. If the divisor $D_r$ is non-empty, then it is ample, and so it is big
and also $C_r \setminus D_r$ is affine. On the other hand, if $D_r$ is the
trivial divisor, then it is clearly not big, and $C_r \setminus D_r$ is a
projective curve and so not affine. 
\end{proof}

For $k \geq 2$, if the $k$-dimensional locally closed strata are affine, then
$\X$ is robust in dimension $k$, but
not conversely, and neither condition is determined by the dual complex.
However, robustness in dimension 2 is determined by its weak tropical complex,
for which we need the following definition.

\begin{defn}[Def.~2.7 in \cite{cartwright-complexes}]
If $q$ is an $(n-2)$-dimensional simplex, the \defi{local intersection matrix}
of~$\Delta$ at~$q$ is the symmetric matrix $M_q$ whose rows and columns are
indexed by the ridges~$r$ containing~$q$, and such that the entry in row~$r$ and
column~$r'$ is:
\begin{equation*}
(M_{q})_{r,r'} = \begin{cases}
\#\{\mbox{facets containing $r$ and $r'$}\} & \mbox{if } r \neq r' \\
-\alpha(v, r) & \mbox{if } r = r',
\end{cases}
\end{equation*}
where $v$ refers to the vertex of~$r$ not contained in~$q$.
A \defi{tropical complex} is a weak tropical complex~$\Delta$ such that the
local intersection matrix $M_q$ has exactly one positive eigenvalue for each
$(n-2)$-dimensional simplex~$q$ of~$\Delta$.
\end{defn}

\begin{prop}[Prop.~2.9 in \cite{cartwright-complexes}]
\label{p:robust-dim-2}
Let $\X$ be a degeneration and $\Delta$ its weak tropical complex. Then $\X$ is
robust in dimension~2 if and only if $\Delta$ is a tropical complex.
\end{prop}

We now prove Proposition~\ref{p:blow-ups}, showing that any projective
degeneration can be modified to one whose locally closed strata are affine, and
thus is robust in all dimensions.

\begin{proof}[Proof of Prop.~\ref{p:blow-ups}]
We fix a component $C_v$ of the special fiber~$\X_0$.
By assumption, $C_v$ is projective, so, by Bertini's theorem, we can
choose smooth and irreducible elements $H_1, \ldots, H_n$ from the linear system
of a very ample divisor on~$\X_0$, such that the $H_i$ intersect both each other
and~$D_v$ transversely. We now blow-up the points of the intersection $H_1 \cap
\cdots \cap H_n$. Then, we blow-up, for each integer $k$ from $1$ to~$n-1$, the
strict transforms of the $k$-dimensional varieties in $H_{i_1} \cap \cdots \cap
H_{i_{n-k}}$ for all indices $1 \leq i_1 <  \ldots < i_{n-k} \leq n$. For each
$k$, these strict transforms are disjoint and thus the order of the blow-ups
within a fixed dimension doesn't matter.

We then repeat the above blow-ups at the strict transform of each component of
$\X$ to obtain the degeneration $\X'$ from the statement. Since $\X'$ is the
iterated blow-up of a regular degeneration at smooth subvarieties, it
is also regular. The centers of the blow-ups intersect the singular locus of the
special fiber transversally, so the special fiber~$\X'_0$ is also reduced and
simple normal crossing. Thus, $\X'$ is a strict degeneration.

In order to show that the locally closed strata of~$\X'$ are all affine, we
first consider the case of the case of a stratum~$C_s'$ of~$\X_0'$ which maps
birationally onto its image in~$\X$. We let $C_s$ be its image in $\X$, and then
$C_s'$ is formed by blowing up $C_s$ at the restrictions of intersections of
very ample divisors from a component~$C_v$ for each vertex~$v$ in~$C_s$. Each of
these blow-ups produces a new component for~$\X'$ and thus the difference $C_s'
\setminus D_s'$ is an open subset of $C_s$ minus the very ample divisors. This
containment may be proper because of the intersection of $C_s$ with other
components of $\X_0$. However, the complement of a very ample divisor is affine,
and the complement of a Cartier divisor in an affine variety is affine, so $C_s'
\setminus D_s'$ is affine, as desired.

Now, we consider the components introduced by the blow-ups. The center for such
a blow-up is a variety~$Y$ within a single component~$C_v$. Since we've blown up
the intersection of $Y$ with very ample divisors in the previous step, the
complement $Y \setminus (Y \cap D_v)$ is affine, as in the previous paragraph.
The blow-up of $Y$ is a projective bundle over $Y$ which intersects $C_v$ along
a section. Let $C_w$ denote this blow-up and then $C_w \setminus D_w$ is an
affine space bundle over an affine variety $Y \setminus (Y \cap D_v)$, and thus
it is affine. Further blow-ups only intersect $C_w$ along $D_w$ and therefore do
not affect the difference $C_w \setminus D_w$, which remains affine in $\X'$. We
conclude that the locally closed strata of~$\X'$ are all affine.
\end{proof}

We now come to divisors on $X$ and their relationship to the combinatorics of
the weak tropical complex $\Delta$. Let $D$ be any divisor on $X$, and let
$\overline D$ denote its closure in $\X$. Since $\X$ is regular, $\overline D$
is a Cartier divisor on $\X$, and we can define the following formal sum of
ridges of $\Delta$:
\begin{equation}\label{eq:specialization}
\rho(D) = \sum_{r \in \Delta_{n-1}} (\degree \overline D \cdot C_r) [r],
\end{equation}
which we call the \defi{coarse specialization} of~$D$. In addition, there is a
refined specialization, which will be construction in
Section~\ref{s:subdivisions}. We now describe the machinery of linear and
piecewise linear functions on $\Delta$, which will be used to define divisors
and linear equivalence.

We first need to generalize beyond just formal sums of the ridges to formal sums
of $(n-1)$-dimensional polyhedra, possibly contained in the interior of a facet.
To be precise, we first identify a single facet~$f$ of~$\Delta$ with the
standard simplex in $\RR^n$, which is the convex hull of the origin together
with the $n$ coordinate vectors. This identification is not unique, but a
permutation of the vertices induces an affine linear automorphism of~$f$, with
integral coefficients, and so it does not affect the places where these
coordinates are used. We define an \defi{$(n-1)$-dimensional polyhedron}
in~$\Delta$ to be an $(n-1)$-dimensional polyhedron of a single facet~$f$
of~$\Delta$, which is defined by linear inequalities with rational coefficients
and real constants, under the above identification.

By a \defi{formal sum of $(n-1)$-dimensional polyhedra}, we will mean a finite,
integral
sum of $(n-1)$-dimensional polyhedra.
From now on, we will drop the dimension from our terminology, because the only
formal sums of polyhedra in this paper will be $(n-1)$-dimensional. Two formal
sums of polyhedra are \defi{equivalent} if they differ by an element of the
subgroup generated by the set of formal sums $[P] - [Q_1] - \cdots [Q_r]$, where
$P$ is any $(n-1)$-dimensional polyhedron, and $Q_1, \cdots, Q_r$ are polyhedra
such that $P = Q_1 \cup \cdots \cup Q_r$ and $Q_i \cap Q_j$ is either empty or a
proper face of both $Q_i$ and $Q_j$.

One can show that any formal sum of polyhedra is equivalent to one
where the terms form the maximal cells of a polyhedral complex. While we will
not use that result in this paper, it in particular means that,
up to equivalence, we can assume that the intersection between any two terms in
a formal sum of polyhedra is a face, and thus has dimension less than $n-1$, as
in part (1) of the following lemma.

\begin{lem}\label{l:subdivision-effective}
Let $Z = \sum a_i [T_i]$ and $Z' = \sum a_i' [T_i']$ be formal sums of
polyhedra which are equivalent to each other on a weak tropical
complex~$\Delta$.
\begin{enumerate}
\item If the intersection $T_i \cap T_j$ of any two terms of $Z$ has dimension
less than $n-1$, and the coefficients $a_i'$ of $Z'$ are positive, then the
coefficients $a_i$ of $Z$ are non-negative.
\item If both sets of coefficients $a_i$ and $a_i'$ are positive, then
$\bigcup P_i = \bigcup P_i'$.
\end{enumerate}
\end{lem}

\begin{proof}
We first assume that the any two terms of $Z$ intersect in dimension less than
$n-1$ and that the coefficients of $Z'$ are positive. Let $p$ be a point of an
arbitrary
term $T_j$ of $Z$. By definition, difference $Z - Z'$ is a finite sum of
expressions of the form $[P] - [Q_1] - \cdots - [Q_r]$, where $P = Q_1 \cup
\ldots \cup Q_r$, and the intersections $Q_i \cap Q_j$ are proper faces of each.
We can assume that $p$ is not contained in any facet of the polyhedra appearing
in these expressions or facet of the $T_i$. Therefore, if $p$ is in the
polyhedron $P$ of one of these expressions, then $p$ is in exactly one of the
corresponding $Q_i$'s. Therefore, equivalence between $Z$ and $Z'$ doesn't
change the sum of the coefficients of the polyhedra containing $p$, meaning
that: \begin{equation}\label{eq:sum-coefficients}
\sum_{P_i \ni p} a_i = \sum_{P_i' \ni p} a_i'.
\end{equation}
Since the coefficients $a_i'$ are positive, the right hand side of
(\ref{eq:sum-coefficients}) is non-negative. By our assumption on $Z$, the left
hand side consists of just the coefficient $a_j$, which is therefore
non-negative.

Second, we assume that the terms of $Z$ and $Z'$ are arbitrary, but the
coefficients $a_i$ and $a_i'$ are all positive. Again, let $p$ be a point in one
term~$T_j$ of~$Z$, such that $p$ is not contained in any of the facets of the
polyhedra appearing in the equivalence between $Z$ and $Z'$, and we again have
the equality (\ref{eq:sum-coefficients}). The left hand side of
(\ref{eq:sum-coefficients}) is positive since $p \in T_j$, so $p$ must also be
contained in some term $T_k'$ of $Z'$. Since we were able to choose $p$ from a
dense open subset of $\bigcup P_i$, and the union $\bigcup P_i'$ is closed, we
have an inclusion $\bigcup P_i \subset \bigcup P_i'$. By symmetry, we have the
reverse inclusion, and thus the desired equality.
\end{proof}

\begin{defn}\label{d:effective-support}
An equivalence class of formal sums of polyhedra is called \defi{effective} if
it contains a formal sum with all coefficients positive.
\end{defn}

The purpose of Lemma~\ref{l:subdivision-effective}(1) is that it gives a criterion
to check whether or not a give formal sum of polyhedra is effective, by
using a representative such that the intersections between the terms all have
dimension less than $n-1$. Lemma~\ref{l:subdivision-effective}(2)
shows that the following definition is independent of the representative chosen.

\begin{defn}
The \defi{support} of
an effective equivalence class of formal sums of polyhedra is the union of
polyhedra in a representative with positive coefficients.
\end{defn}

For example, if we take $P$ to be an $(n-1)$-dimensional polyhedron, and $Q$ any
$(n-1)$-dimensional polyhedron properly contained in~$P$, then $[P]-[Q]$ is
effective, even though it is not written with positive coefficients, and its
support is $P \setminus Q^o$, where $Q^o$ denotes the relative interior of $Q$.
We can find a subdivision $P = Q \cup Q_2 \cup \cdots \cup Q_r$, as in the
equivalence relation, and then $[P] - [Q]$ is equivalent to $[Q_2] + \cdots +
[Q_r]$, which has positive coefficients, thus showing that $[P]-[Q]$ is
effective, and its support is $P \setminus Q^o$, as claimed.

We now define piecewise linear and linear functions. First, these will allow us
to define divisors as equivalence classes of formal sums of polyhedra satisfying
a balancing condition. Second, piecewise linear functions will be used to define
linear equivalence between divisors, and thus complete linear series.

\begin{defn}[Def. 3.1 in \cite{cartwright-complexes}]\label{d:pl}
A \defi{piecewise linear function} or \defi{PL function} on a weak tropical
complex $\Delta$ is a continuous function $\phi$ such that, restricted to each
simplex of $\Delta$, $\phi$ is a piecewise linear function with integral slopes,
under the identification of the simplex with a standard simplex in $\RR^k$.
\end{defn}

\begin{defn}[Const. 3.2 in \cite{cartwright-complexes}]\label{d:linear}
Let $\Delta$ be an $n$-dimensional weak tropical complex. For any ridge $r$ of
$\Delta$, we define $N_r$ to be the simplicial complex which consists of a
central $(n-1)$-dimensional simplex, together with one $n$-dimensional simplex
for each facet containing $r$, and we let $d$ denote the number of such facdets.
Thus, if $N_r^o$ is the union of the interiors of the $n$-dimensional simplices
and the central $(n-1)$-dimensional simplices, then $N_r^o$ is homeomorphic to a
neighborhood of the interior of $r$.

We can embed $N_r$ in the real vector space $\RR^{d+n}$ by sending $w_i$ to the
$i$th coordinate vector of $\RR^{d+n}$ and $v_i$ to the $(d+i)$th coordinate
vector, where $v_1, \ldots, v_n$ are the vertices of the central
$(n-1)$-dimensional simplex, and $w_1, \ldots, w_d$ are the other vertices of
$N_r$. Then, we define $L_r$ to be the quotient of $\RR^{d+n}$ by the
one-dimensional vector space generated by $(1, \ldots, 1, -\alpha(v_1, r),
\ldots, -\alpha(v_n,r))$. We define $\phi_r \colon N_r \rightarrow L_r$ to be
the composition of the inclusion followed by the projection map.

A \defi{linear} function $\phi$ on an open subset $U$ of a weak tropical complex
$\Delta$ is a PL function such that:
\begin{enumerate}
\item For any facet $f$ of $\Delta$, the restriction $\phi \vert_{U \cap f}$
is the restriction of an affine linear function on $\RR^d$, under the
indentification of $f$ with the standard unit simplex in $\RR^d$.
\item For any ridge $r$ of $\Delta$, the restriction $\phi \vert_{U \cap N_r^o}$
is a composition $\ell \circ \phi_r$, where $\ell \colon L_r \rightarrow \RR$ is
a linear function.
\end{enumerate}
\end{defn}

\begin{defn}[Const. 3.5 and Def. 3.7 in \cite{cartwright-complexes}]\label{d:mult-balancing}
Let $Z$ be a formal sum of $(n-1)$-dimensional polyhedra. Let $Q$ be an
$(n-2)$-dimensional face of a term of $Z$ and assume that every for term $P_i$
of $Z$ which intersects $Z$ contains it as a face. Now let $\phi$ be a linear
function on a neighborhood of the interior of $Q$, which is constant on~$Q$.

For each term $P_i$, use the coordinates in $\RR^n$ coming from the
identification of the facet $f \subset \Delta$ containing $P_i$ with the
standard unit simplex. We let $\mathbf w_i \in \ZZ^n$ be a vector such that
$\mathbf w_i \cdot x$ is constant for all $x$ in $P_i$. We assume that the
entries of $\mathbf w_i$ are relatively prime, which uniquely determines
$\mathbf w_i$ up to a sign. Then let $\mathbf v_i \in \ZZ$ be the normal vector
of a supporting hyperplane of $Q$, meaning that for some constant $c_i \in \RR$,
$\mathbf v_i \cdot x \geq c_i$ for all $x \in P_i$, and $Q$ is the set of $x$
such that $\mathbf v_i \cdot x = c_i$. In addition, we can assume that for each
integer $k$, the entries of $\mathbf v_i + k \mathbf w_i$ are relatively prime,
which determines $\mathbf v_i$ up to a multiple of $\mathbf w_i$. Since $\phi$
is constant on $Q$, $\phi \vert_{P_i}$ can be written as $\phi\vert_{P_i} (x) =
a_i \mathbf v_i \cdot x + t$, with $t \in \RR$ and $a_i$ an integer by the
assumption on $\mathbf v_i$ and the fact that $\phi$ has linear slopes. We refer
$a_i$ as the \defi{slope of $\phi$ along $P_i$}. If $m_i$ is the multiplicity of
$P_i$ in $Z$, then we define the \defi{multiplicity of $\phi$ along $Q$ of $Z$}
to be:
\begin{equation*}
\mult_{Z,Q}(\phi) = \sum a_i m_i.
\end{equation*}

We say that a formal sum of polyhedra~$Z$ is \defi{balanced} if each
$(n-2)$-dimensional face $Q$ of a term of $Z$ is a face of all the terms which
contain it, and $\mult_{Z,Q}(\phi) = 0$ for any linear function $\phi$ an a
neighborhood of the interior of any $(n-2)$-dimensional face $Q$ of a term of
$Z$, which is constant on~$Q$.
\end{defn}

\begin{lem}[Lem. 3.8 in \cite{cartwright-complexes}]\label{l:equiv-balanced}
Let $Z$ and $Z'$ be equivalent formal sums of $(n-1)$-dimensional polyhedra such
that any $(n-2)$-dimensional face of a term of $Z$ is a face of all terms of $Z$
which contain it, and similarly with $Z'$. Then, $Z$ is balanced if and only if
$Z'$ is balanced.
\end{lem}

\begin{defn}[Def. 3.9 in \cite{cartwright-complexes}]\label{d:divisor}
A \defi{(Weil) divisor} on a weak tropical complex is an equivalence class of
formal sums of $(n-1)$-dimensional polyhedra such that some 
representative is balanced.
\end{defn}

In~\cite{cartwright-complexes}, Cartier divisors on weak tropical complexes are
also defined, which are the appropriate tool for intersection theory. With a
exception of the Riemann-Roch conjecture at the end of
Section~\ref{s:linear-series}, only Weil divisors are relevant in this paper, which we
usually refer to as \defi{divisors}. In particular, divisors arise from the
specialization of divisors:

\begin{prop}[Prop. 3.16 in \cite{cartwright-complexes}]\label{p:specialization-weil}
Let $\X$ be a degeneration which is robust in dimension $2$ and let $\Delta$ be
its tropical complex. If $D$ is a divisor on $X$, then $\rho(D)$ is a divisor on
$\Delta$.
\end{prop}

Linear equivalence between divisors on a weak tropical complex is analogous to
linear equivalence on an algebraic variety, with PL functions taking the place
of rational functions. The following proposition characterizes the divisor
associated to a PL function, which closely links them to the linear functions
from Definition~\ref{d:linear}.

\begin{prop}[Prop. 4.1 in \cite{cartwright-complexes}]\label{p:div}
Let $\Delta$ be a weak tropical complex and assume that every ridge is contained
in a facet. There is a unique function $\div$ from PL functions on  an open
subset $U \subset \Delta$ to equivalence classes of formal sums of polyhedra on
$U$ such that:
\begin{itemize}
\item[(i)] For any PL functions $\phi$ and $\phi'$ on $U$, $\div(\phi+\phi') =
\div(\phi) + \div(\phi')$.
\item[(ii)] If $V \subset U$ are open sets, and $\phi$ is a PL function on $U$,
then $\div(\phi \vert_V) = \div(\phi)\vert_V$.
\item[(iii)] A function $\phi$ is linear if and only if $\div(\phi)$ is trivial.
\item[(iv)] If $\phi$ is identically zero outside of a single facet $f$ of
$\Delta$, and $\phi$ is defined as:
\begin{equation*}
\phi(x) = \max\{ \lambda \cdot x, 0\},
\end{equation*}
where $x$ is a coordinate vector, using the identification of $f$ with a
standard unit simplex in $\RR^n$, and $\lambda$ is an integral vector whose
entries have no common divisor, then $\div(\phi) = [H]$, where $H = \{x \in f
\cap U \mid \lambda \cdot x = 0\}$, using the same coordinates.
\end{itemize}
\end{prop}

In \cite{cartwright-complexes}, Proposition 4.1 also has an additional
normalization assumption in order to apply to ridges that aren't contained in a
facet. However, the degenerations of interest in this paper are robust in
dimension 1, and thus, by Proposition~\ref{p:robust-dim-1}, every ridge is
contained in a facet.

\begin{defn}[Def. 5.1 in \cite{cartwright-complexes}]\label{d:equivalence}
Two divisors $D$ and $D'$ on a weak tropical complex~$\Delta$ are \defi{linearly
equivalent} if $D-D' = \div(\phi)$ for some PL function $\phi$ on $\Delta$.
\end{defn}

\begin{prop}[Prop. 5.2 in \cite{cartwright-complexes}]\label{p:lin-equivalent}
Let $\X$ be a degeneration which is robust in dimension $2$, with $\Delta$ as
its tropical complex. If $D$ and $D'$ are linearly equivalent divisors on $X$,
then $\rho(D)$ is linearly equivalent to $\rho(D')$.
\end{prop}

\section{Linear series}\label{s:linear-series}

A \defi{complete linear series} of a divisor~$D$ is the set of all effective
divisors linearly equivalent to~$D$. The invariant $h^0(\Delta, D)$ appearing in
Theorem~\ref{t:specialization} is a measure of how large the complete linear
series of~$D$ is, analogous to the dimension of the global sections of $\mathcal
O(D)$ on an algebraic variety. In this section, we define this invariant and
give some examples and applications.

\begin{defn}\label{d:h0}
We say that a point $p \in \Delta$ is \defi{rational} if its coordinates are
rational when we identify a $k$-simplex containing $p$ with a standard unit
simplex in~$\RR^k$. We define $h^0(\Delta, D)$ to be the cardinality~$m$ of the
smallest set of rational points $p_1, \ldots, p_m$ such that there is no
effective divisor $D'$ linearly equivalent to $D$ such that $D'$ contains $p_1,
\ldots, p_m$. If there is no such finite set of points, then $h^0(\Delta, D)$ is
defined to be $\infty$.
\end{defn}

The rationality condition on the points in Definition~\ref{d:h0} is needed
for technical reasons in the proof of Theorem~\ref{t:specialization} and we
expect that it can be dropped without changing the definition.

\begin{conj}\label{conj:non-rational-points}
The definition of $h^0$ is equivalent to the analogous definition where the
$p_i$ are allowed to be arbitrary, not necessarily rational, points.
\end{conj}

In the case of a $1$-dimensional tropical complex~$\Gamma$, the quantity
$h^0(\Gamma, D)$ is essentially the same as the rank of the divisor as
introduced by Baker and Norine~\cite{baker-norine}, and extended to metric
graphs in~\cite{gathmann-kerber} and~\cite{mikhalkin-zharkov}, with the
exception that our convention differs from theirs by~$1$. While our definition
requires the points to be distinct and rational, these restrictions don't affect
the definition, by the following case of
Conjecture~\ref{conj:non-rational-points}:

\begin{prop}\label{p:agree-curve}
If $\Gamma$ is a $1$-dimensional tropical complex with at least one edge and
$r(D)$ is the rank of a divisor~$D$ as in~\cite{baker-norine}, then $h^0(D) =
r(D) + 1$.
\end{prop}

\begin{proof}
Recall that the rank of $D$ is the largest number $r$ such that for any $r$
points $p_1, \ldots, p_r$ in $\Gamma$, the difference $D - [p_1] - \ldots -
[p_r]$
is linearly equivalent to an effective divisor~$D'$. Thus, $D' + [p_1] + \ldots
+ [p_r]$ is an effective divisor linearly equivalent to $D$ and it clearly
contains
the $p_i$, so $h^0(D) \geq r(D) + 1$.
To show the reverse inequality, we assume that
 there exist points $p_1, \ldots, p_{r(D) + 1}$ such
that $D$ is not linearly equivalent $D' + [p_1]
+ \cdots + [p_{r(D) + 1}]$ for any effective divisor $D'$. In other words, if we
write $\lvert D \rvert$ for the subset of~$\Gamma^d$ which are effective
divisors linearly equivalent to $D$,
then $\lvert D \rvert$ does not intersect the subset
\begin{equation*}
\{p_1\} \times \cdots \times \{p_{r(D) + 1}\} \times \Gamma \times \cdots \times
\Gamma \subset \Gamma^d.
\end{equation*}
However, $\lvert D \rvert$ is a closed
subset~\cite{mikhalkin-zharkov}*{Thm.~6.2}. Therefore, since $\Gamma$ is not a
point, we can perturb the points $p_i$ slightly and the intersection with
$\lvert D \rvert$ will still be empty. In particular, we can make the $p_i$
distinct and rational. Thus, $h^0(D) \leq r(D) + 1$, so we've proved
the proposition.
\end{proof}

\begin{figure}
\includegraphics{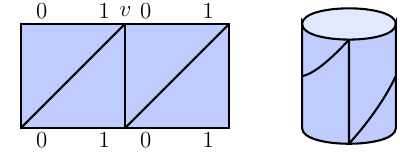}
\caption{The tropical complex~$\Delta$ used in Theorem~\ref{t:no-lift} is
formed by identifying the edges on the left and
right from the triangulation on the left to form the cylinder shown on the right.
For the edges forming the top and bottom circles, the structure constants are
indicated by the numbers adjacent to each endpoint. The other edges all have
degree~$2$ and all structure constants on these edges are taken to
be~$1$.}\label{f:hopf-regular}
\end{figure}

We now give some applications of Theorem~\ref{t:specialization} for
2-dimensional tropical complexes, beginning with the proof of
Theorem~\ref{t:no-lift}, that there exists a 2-dimensional tropical complex
which doesn't lift. For our example, we take the cylinder depicted in
Figure~\ref{f:hopf-regular}, which is a variant of Example~4.6
in~\cite{cartwright-surfaces}, with the circumference of the cylinder increased
to~$2$ so that the underlying complex is a regular $\Delta$-complex. This
tropical complex is a tropical analogue of Hopf surface, which are non-algebraic
complex surfaces. Tropical Hopf manifolds are explored at greater length, and in
a different framework, in~\cite{ruggiero-shaw}.

\begin{lem}\label{l:series-cylinder}
Let $\Delta$ be the tropical complex in Figure~\ref{f:hopf-regular} and let $D$
be the sum of the top two edges. Then any divisor linearly equivalent to $mD$ is
the sum of $m$ copies of parallel circles, and thus $h^0(\Delta, mD) = m+1$.
\end{lem}

\begin{proof}
Let $\phi$ be a PL function on~$\Delta$ such that $\div(\phi) + mD$ is
effective. Let $C$ be a horizontal circle on the cylinder~$\Delta$ parallel
to~$D$, but not equal to~$D$ or the bottom of the cylinder. Let $p$ be a point
of~$C$ at which the restriction $\phi \vert_C$ achieves its maximum. The key
point is that the linear functions on a neighborhood~$U$ of~$p$ embed $U$ in an
open subset of $\RR^2$ and $C$ remains a straight line in this embedding. Since
the divisor of $\phi$ is effective in $U$, $\phi\vert_{U}$ is a convex function
of~$\RR^2$. Thus, the restriction to $C$, which is a line segment, is also
convex, but $\phi \vert_C$ achieves its maximum at~$p$ and so $\phi \vert_{C
\cap U}$ must be constant. Since $\phi$ is constant in any neighborhood of a
point in~$C$ where $\phi\vert_C$ achieves its maximum, $\phi$ is constant
on~$C$.

Thus, $\phi$ is a function solely of the vertical coordinate on~$\Delta$, and
$\phi$ defines a linear equivalence between $mD$ and $m$ horizontal circles,
possibly not distinct. For any $m$ points, these circles can be chosen to
contain them, but for $m+1$ points chosen to have $m+1$ distinct heights, no sum
of $m$ circles will contain all of the points. Therefore, $h^0(\Delta, mD) =
m+1$, as claimed.
\end{proof}

\begin{proof}[Proof of Theorem~\ref{t:no-lift}]
Suppose that $\X$ is a degeneration over a discrete valuation ring whose weak
tropical complex is~$\Delta$. Since all the maximal simplices of $\Delta$ are
$2$-dimensional, and it can be checked that $\Delta$ is a tropical complex, $\X$
must be robust in dimensions $1$ and~$2$ by Propositions~\ref{p:robust-dim-1}
and~\ref{p:robust-dim-2}, respectively. The general fiber $X$ is a smooth proper
surface and therefore projective. We let $A$ be an ample divisor on $X$. By
Proposition~\ref{p:specialization-weil}, $\rho(A)$ is a divisor, and from the
definition of the specialization~(\ref{eq:specialization}), we know that
$\rho(A)$ is a linear combination of the ridges of $\Delta$.

We now examine the possibilities for $\rho(A)$.
The set of all formal sums of the 8 ridges of $\Delta$ forms a free Abelian
group of rank 8, and those satisfying the
balancing condition are a subgroup. In a neighborhood of each vertex $v$ of
$\Delta$, the space of linear functions can be computed by considering functions
which are linear on the interior of each 2-dimensional simplex, and satisfy the
condition in Definition~\ref{d:linear} along each ridge containing $v$. In each
case,
the space of linear functions have only one parameter for the slopes.
For example, let $v$ be the vertex indicated in Figure~\ref{f:hopf-regular}, and
$e_1$, $e_2$, $e_3$, and $e_4$ be the edges incident to $v$, beginning with
$e_1$ to the left and continuing counterclockwise. Then, the linear functions
are constant along the diagonal edge $e_2$, and affine linear under the
embedding shown.
Therefore, the slopes along $e_3$ and $e_4$ are equal, and are
the negative of the slope along $e_1$. Therefore, if $c_1$, $c_3$, and $c_4$ are
the coefficients of these three edges in a divisor, then $c_1 - c_3 -c_4 = 0$.
Similarly,
the other 3 vertices of $\Delta$ also give linear conditions on the coefficients
of a divisor, and one can check that these are independent, so the group of
formal sums of the ridges which are divisors is a free Abelian group of rank 4.

Next, one can compute the subgroup of divisors which are linearly equivalent to the
trivial one by considering PL functions which are linear on the interior of each
$2$-dimensional simplex. Using the properties in Proposition~\ref{p:div}, one
can then check that the group of divisor classes which are linearly equivalent
to a formal sum of ridges is isomorphic to $\ZZ$, with the divisor $D$ from
Lemma~\ref{l:series-cylinder} representing twice a generator.

Therefore, $\rho(2A) = 2 \rho(A)$ is linearly equivalent to $lD$ for some
integer $l$. Thus, by Theorem~\ref{t:specialization}, and
Lemma~\ref{l:series-cylinder} we have the inequality
\begin{equation*}
h^0(X, \mathcal O(2mA)) \leq h^0(\Delta, mlD) = ml + 1.
\end{equation*}
In particular, $h^0(X, \mathcal O(2mA))$ grows at most linearly in $m$, which
contradicts the fact that $A$ is ample divisor on a surface, and so its global
sections grow quadratically. Therefore, we conclude that no degeneration
$\X$ with weak tropical complex $\Delta$ can exist.
\end{proof}

The following example shows one case where the inequality in
Theorem~\ref{t:specialization} can be sharp for a divisor and all of its
multiples.

\begin{ex}\label{ex:k3-series}
We consider a degeneration of a quartic surface which will yield a tetrahedron
as its dual complex.
We start with the variety $\widetilde \X
\subset \PP^3_R$, where $R = \CC[[t]]$, which is defined by the determinant:
\begin{equation}\label{eq:k3-series}
\begin{vmatrix}
xy & (3x^2 + 2 y^2 + z^2 - w^2) t \\
3x^2 + y^2 + 2z^2 + w^2 &
zw + t(x^2 + y^2 + z^2 + w^2)
\end{vmatrix}
\end{equation}
The reason for this determinantal form is that we can immediately see that the
general fiber $X$ contains the scheme defined by the equations in
the bottom row of the matrix~(\ref{eq:k3-series}), and one can check that this
is a smooth complete intersection and therefore an elliptic curve~$E$. The
curve~$E$ is part of a pencil interpolating between $E$ and the curve defined by
the equations in the top row of~(\ref{eq:k3-series}). Thus, the complete linear
series of~$E$ defines a map to $\PP^1_K$ with connected fibers, and so
$h^0(\widetilde \X_K, \mathcal O(mE)) = m+1$ for all $m \geq 0$.

Now we want to show how Theorem~\ref{t:specialization} gives a sharp
bound for this value of $h^0(\widetilde \X_K, \mathcal O(mE))$. We can rewrite
(\ref{eq:k3-series}) as $xyzw + tf$, where
\begin{equation*}
f = xy(x^2+y^2+z^2+w^2) - (3x^2+2y^2+z^2-w^2)(3x^2+y^2+2z^2+w^2),
\end{equation*}
to see that the special fiber of $\widetilde X$ is the union of the four
coordinate planes. However, $\widetilde X$ is not a regular degeneration, but
has 24 ordinary double point singularities at the common intersection of the
quartic~$f$ and the 6 coordinate lines in $\PP^3$. We wish to resolve each
singularity without introducing any new components in the special fibers,
which can be done by blowing up one of the two planes containing the
singularity.

As in \cite{cartwright-complexes}*{Ex.~2.10}, we can obtain a
symmetric tropical complex by blowing up one plane at 2 of the 4 singularities
along each line and blowing up the other plane at the other 2 singularites, but
we pay special attention to the four lines defined by the intersection of one of
the planes $x=0$ or $y=0$ with either $z=0$ or $w=0$. Note that when either $x$
or $y$ vanishes, $f$ factors as the product of two quadrics, which we write as
as
\begin{equation*}
g = 3x^2 + 2y^2 + z^2 - w^2
\quad\mbox{and}\quad h = 3x^2 + y^2 + 2z^2 +w^2.
\end{equation*}
At each point
of intersection of one of these four lines with $g=0$, we blow up either the
$x=0$ or $y=0$, as appropriate, and at each point of intersection with $h=0$, we
blow up either the $z = 0$ or $w=0$ plane. Along the other 2 coordinate lines,
we blow up each containing plane at 2 of the singularities, but chosen
arbitrarily.

After these blow ups, the resulting scheme $\X$ is a regular semistable
degeneration. The special fiber has four components, each of which is the
blow up of $\PP^2$ at 6 points, 2 on each coordinate line. The 1-dimensional
strata in each component are the union of the strict transforms of the
coordinate lines, and one can check that the union of these lines is a big
divisor, so $\X$ is robust and thus we can apply the specialization inequality
from Theorem~\ref{t:specialization}. The tropical complex~$\Delta$ is a
tetrahedron with all structure constants $\alpha(v,e)$ equal to $1$, using the
fact that the strict transform of each coordinate line has self-intersection
$-1$.

We first calculate the specialization $\rho([E])$, where $E$ is the elliptic
curve on
$X$ defined above. In the original singular degeneration~$\widetilde\X$, the
closure~$\overline E$ of~$E$ degenerates to the union of two conics, in the
$z=0$ and $w=0$ planes, respectively, each defined by the restriction of the
polynomial~$h$. However, we chose to blow up the $z=0$ and $w=0$ planes at the
intersections of $h$ with the $x=0$ and $y=0$ lines, which removes the
intersections between $\overline E$ and those lines. Therefore, the
only 1-dimensional stratum of $\X$ that $\overline E$ intersects is the $x=y=0$
line, which it does with multiplicity 2. Thus, $D = \rho(E)$ is twice one edge
of the tetrahedron~$\Delta$.

Finally, to justify our claim that the specialization inequality is sharp, we
need to show that $h^0(\Delta, mD) = m+1$ for all non-negative integers $m$. For
this, we start with a PL function~$\phi$ which is constant on the support of~$D$
and increases with slope 1 on each of the two containing facets, which
establishes a linear equivalence between $D$ and any cycle on the tetrahedron
parallel to~$D$. Moreover, if we ignore both the edge supporting~$D$ and the
edge opposite it in~$\Delta$, we are left with an open cylinder, which has the
same theory of linear equivalence as if we removed the top and bottom circles
from the cylinder in Figure~\ref{f:hopf-regular}, used in the proof of
Theorem~\ref{t:no-lift}. Therefore, Lemma~\ref{l:series-cylinder} shows any
divisor in the complete linear series of $mD$ consists of $kD + k'D'$ plus $l$
cycles parallel to~$D$, where where $D'$ is twice the edge opposite~$D$ and $k +
k' + l = m$. As in the proof of Lemma~\ref{l:series-cylinder}, such divisors
can be chosen to contain any $m$ points, but not $m+1$
points at different distances from~$D$, and so $h^0(\Delta, mD) = m+1$.
\end{ex}

The sharpness of the inequality in Example~\ref{ex:k3-series} depended on our
choice of blow-ups in constructing the resolution~$\X$. For surfaces and
higher-dimensional varieties, unlike the case of curves, there is no single
minimal regular semistable degeneration, and the choice of the model can affect
the bounds from Theorem~\ref{t:specialization}.

We close this section by stating a conjectural form of Riemann-Roch for
2-dimensional tropical complexes. Although we have no definition for higher
cohomology, we can assume Serre duality to justify taking $h^0(\Delta, K_\Delta
- D)$ as a replacement for the top cohomology. Here, we take $K_\Delta$ of any
tropical complex~$\Delta$ to be the formal sum of polyhedra
\begin{equation*}
K_\Delta = \sum_{r} (\degree r - 2)[r],
\end{equation*}
where $r$ ranges over
the ridges~$r$ of~$\Delta$, generalizing the definition for curves. In
dimension~1, this approach yielded the Riemann-Roch theorems for graphs and
tropical curves~\cites{baker-norine,gathmann-kerber,mikhalkin-zharkov}.

In addition, we need two concepts on tropical complexes which we were not
defined in Section~\ref{s:tropical-complexes} and are not needed in the rest of
the paper. First, \cite{cartwright-complexes}*{Def.~4.2} defines Cartier
divisors on a weak tropical complex, and every Cartier divisor on a
2-dimensional weak tropical complex is a Weil
divisor~\cite{cartwright-surfaces}*{Prop.~2.8}. Second, 2-dimensional tropical
complexes have a symmetric, bilinear pairing on divisors, which is integral on
Cartier divisors~\cite{cartwright-complexes}*{Prop.~4.7}
and~\cite{cartwright-surfaces}*{Prop.~2.9}. With these ingredients, we
conjecture the following, which is a formal analogue to the Riemann-Roch theorem
for algebraic surfaces, as it was stated before the introduction of cohomology:

\begin{conj}[Riemann-Roch]\label{conj:riemann-roch}
Let $D$ be a Cartier divisor on a 2-dimensional tropical complex~$\Delta$ and
assume that $K_\Delta$ is also a Cartier divisor. Then
\begin{equation}\label{eq:r-r}
h^0(\Delta, D) + h^0(\Delta, K_\Delta - D) \geq
\frac{\degree D \cdot (D - K_\Delta)}{2} + \chi(\Delta),
\end{equation}
where $\chi(\Delta)$ is the Euler characteristic of underlying topological space
of~$\Delta$.
\end{conj}

\begin{ex}\label{ex:tetra-rr}
We verify Conjecture~\ref{conj:riemann-roch} for the tetrahedron~$\Delta$
from Example~\ref{ex:k3-series} and for any multiple $mD$,
where $D$ is twice an edge of the tetrahedron. Since every edge is contained in
exactly two facets, $K_\Delta$ is trivial. Moreover, Example~\ref{ex:k3-series}
showed that $D$ is linearly equivalent to a parallel circle, which is disjoint
from $D$, and so $\deg D^2 = 0$.
Therefore, the right hand side of~(\ref{eq:r-r}) is
$\chi(\Delta) = 2$ in this example.

Example~\ref{ex:k3-series} computed that $h^0(\Delta, mD) = m+1$ for $m \geq 0$.
The same method as in that example shows that if $m > 0$, then $h^0(\Delta, K_\Delta
-mD) = h^0(\Delta, -mD) = 0$. (That the negative of an effective
divisor is non-effective also holds for a broad class of 2-dimensional tropical
complexes, as a consequence of \cite{cartwright-surfaces}*{Cor. 2.12}.) Thus,
for any integer $m$, the left hand side of~(\ref{eq:r-r}) is:
\begin{equation*}
h^0(\Delta, mD) + h^0(\Delta, -mD) = \begin{cases}
m + 1 & \mbox{if } m \geq 1 \\
2 & \mbox{if } m = 0 \\
-m + 1 & \mbox{if } m \leq -1,
\end{cases}
\end{equation*}
and so the Riemann-Roch inequality is satisfied.

Note that $\frac{1}{2}D$ still has integral coefficients and so is a
Weil divisor. Moreover, one can check that $h^0(\Delta, \frac{1}{2} D)
= 1$ and $h^0(\Delta, -\frac{1}{2}D) = 0$ and thus the Riemann-Roch inequality
is not satisfied for $D$, which doesn't contradict
Conjecture~\ref{conj:riemann-roch} because $\frac{1}{2} D$ is not Cartier.
\end{ex}

\section{Subdivisions}\label{s:subdivisions}

In this section, we study subdivisions of a weak tropical complex, which
correspond to ramified extensions of the discrete valuation ring~$R$ followed by
toroidal resolution of singularities. Here, the key technology is the
relationship between toroidal maps and polyhedral subdivisions as
in~\cite{kkmsd}. Combinatorially, passing to a degree~$m$ totally ramified
extension of the DVR~$R$ corresponds to scaling the simplices of~$\Delta$ by a
factor of~$m$, and the toroidal resolution is given by a unimodular subdivision
of the scaled simplices. The importance of subdivisions for us is as a tool when
we define the refined specialization in the next section.

\begin{const}[Subdivisions]\label{const:subdivision}
By an \defi{order-$m$ subdivision} of a weak tropical complex~$\Delta$, we mean
a $\Delta$-complex which is formed by replacing each simplex of~$\Delta$ by a
unimodular integral subdivision of the standard simplex scaled by~$m$, together
with the structure constants obtained as follows.

First suppose that $r'$ is a ridge of $\Delta'$ meeting the interior of a
facet of~$\Delta$.
In particular, if $v_1', \ldots, v_n'$
are the vertices of $r'$, and $w_1$ and $w_2$ contained in the two facets
containing $r'$, then the midpoint $(w_1 + w_2) / 2$ is contained in the plane
spanned by $r'$. Thus, we can write:
\begin{equation*}
(w_1 + w_2) / 2 = c_1 v_1' + \cdots + c_n v_n'
\end{equation*}
for some coefficients $c_1, \ldots c_n$ with $c_1 + \cdots + c_n = 1$. We set
$\alpha(v_i', r') = 2 c_i$.

On the other hand, suppose that $r'$ is a ridge of~$\Delta'$ which is contained
in a ridge~$r$ of~$\Delta$, and $r$ has degree~$d$. We represent the points of
each facet of~$\Delta$ which contains $r$ by $(n+1)$ non-negative coordinates
with total sum equal to~$m$. In the $i$th such facet, the unique facet of the
subdivision~$\Delta'$ containing $r'$ is the convex hull of~$r'$ and a unique
point, whose coordinate is $(x_{i,1}, \ldots, x_{i,n}, 1)$, by unimodularity.
Likewise, the points of~$r$ can be given coordinates consisting of $n$
non-negative  real numbers, also summing to~$m$. In such coordinates, we
represent the $i$th vertex~$v_i$ of $r'$ by the vector $(y_{i,1}, \ldots,
y_{i,n})$. Finally, if $v_1, \ldots, v_n$ are the vertices of~$r$, we determine
the structure constants of~$r'$ by the equation:
\begin{equation}\label{eq:subdivision}
\begin{pmatrix}
\alpha(v_1', r') \\ \vdots \\ \alpha(v_n', r')
\end{pmatrix}
= \begin{pmatrix}
y_{1,1} & \cdots & y_{n, 1} \\
\vdots & & \vdots \\
y_{1,n} & \cdots & y_{n,n}
\end{pmatrix}^{-1} \!
\begin{pmatrix}
\alpha(v_1, r) + x_{1,1} + \cdots + x_{d,1} \\
\vdots \\
\alpha(v_n, r) + x_{1, n} + \cdots + x_{d,n}
\end{pmatrix}\!.
\end{equation}
That these structure constants form a weak tropical complex is
verified in Proposition~\ref{p:subdivision}, after the following example.
\end{const}

\begin{ex}
\begin{figure}
\includegraphics{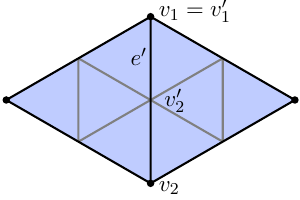}
\caption{The gray lines show a order~$2$ subdivision of two triangles glued
along a common edge. In both the original and subdivided complex, $\alpha(v,r) = 1$
for all ridges~$r$ which are contained in two facets.}\label{f:subdivision}
\end{figure}
Consider the order~$2$ subdivision of the complex~$\Delta$ shown in
Figure~\ref{f:subdivision}.
Let $e'$ be the upper segment of the central edge, with $v_1'$ the top vertex
and $v_2'$ the middle vertex of the edge. If $v_1$ is also the top vertex and
$v_2$ is the bottom vertex, then $(y_{1,1}, y_{1,2}) = (2,0)$ and $(y_{2,1},
y_{2,2}) = (1,1)$. The vertices adjacent to $e'$ both have coordinate $(1, 0,
1)$. Thus, applying~(\ref{eq:subdivision}), we get:
\begin{equation*}
\begin{pmatrix} \alpha(v_1', e') \\ \alpha(v_2', e') \end{pmatrix} =
\begin{pmatrix}
2 & 1 \\ 0 & 1
\end{pmatrix}^{-1}
\begin{pmatrix}
1 + 1 + 1 \\ 1 + 0 + 0
\end{pmatrix} 
=
\begin{pmatrix}
{\textstyle\frac{1}{2}} & {\textstyle - \frac{1}{2}} \\ 0 & 1
\end{pmatrix}
\begin{pmatrix}
3 \\ 1
\end{pmatrix}
= \begin{pmatrix} 1 \\ 1 \end{pmatrix}
\end{equation*}
More generally, we have $\alpha(v,e) = 1$ for every internal edge in the
subdivided complex~$\Delta'$.
\end{ex}

\begin{prop}\label{p:subdivision}
The subdivision in Construction~\ref{const:subdivision} results in a weak
tropical complex.
\end{prop}

\begin{proof}
We need to show that the $\alpha(v', r')$ from
Construction~\ref{const:subdivision} are integral and that they satisfy the
required identity~(\ref{eq:ridge-identity}) for each ridge.
We retain the same notation as in Construction~\ref{const:subdivision} and
we first consider a
ridge $r'$ meeting the interior of a facet of~$\Delta$. Then the coordinates of
$(w_1 + w_2)/2$ will be half-integers, and so the $c_i$ are also half-integers,
and $\alpha(v_i', r') = 2c_i$ is an integer. Moreover, the sum of the
$\alpha(v_i', r')$ will equal $2$, which is the degree of~$r'$ in $\Delta'$.

Now consider a ridge $r'$ of the subdivision~$\Delta'$ contained in a ridge~$r$
of~$\Delta$. The fact that the vectors $v_1', \ldots, v_n'$ form the vertices of
one simplex of a unimodular triangulation imply that the differences between
pairs of the vectors $(y_{i,1}, \ldots, y_{i,n})$ span the integral vectors
$\ZZ^{n}$ whose sum is $0$. Moreover, each vector satisfies the affine linear
equation $y_{i,1} + \cdots + y_{i,n} = m$, so the span of these vectors includes
all vectors in $\ZZ^{n}$ whose sum is a multiple of~$m$. Thus, to check that the
coefficients defined in~(\ref{eq:subdivision}) are integral, it is sufficient to
check that the sum of the entries of the vector from the right hand side of
(\ref{eq:subdivision}) is a multiple
of~$m$. Indeed, we know that $\alpha(v_1, r) + \cdots + \alpha(v_n, r) = d$ and
that $x_{i, 1} + \cdots + x_{i,n} = m-1$, so the total sum is $dm$. Moreover,
this implies that $\alpha(v_1', r') + \cdots \alpha(v_n', r') = d$, which
completes the proof that the subdivision is a weak tropical complex.
\end{proof}

It is immediate from the definition that a subdivision of a weak tropical
complex is homeomorphic to the original complex. In addition, we now show that
the subdivision behaves almost identically with the properties of weak tropical
complexes introduced in Section~\ref{s:tropical-complexes}, namely linear
functions and divisors.

\begin{lem}\label{l:subdivide-linear}
Let $\Delta'$ be a order~$m$ subdivision of a weak tropical complex~$\Delta$,
and we identify the realizations of~$\Delta'$ and~$\Delta$ via the natural
homeomorphism. Then a function~$\phi$ on $\Delta$ is linear (respectively PL) if
and only if $m\phi$ is linear (respectively PL) on $\Delta'$.
\end{lem}

\begin{proof}
The agreement of PL functions is clear, so long as the scaled function $m\phi$
on $\Delta'$ has integral slopes on its domains of linearity. However, each
facet $f'$ of $\Delta'$ is identified with a $(1/m)$-scaled copy of the standard
unit simplex. If we instead identify $f'$ with the standard unit simplex, then
the slopes of $\phi$ will be scaled by $1/m$, and so the slopes of $m\phi$ will
all be integers.

For linear functions, we let $r'$ be a ridge of of $\Delta'$ and we need to show
that linear functions on a neighborhood of $r'$ in $\Delta$, scaled by $m$,
agree with linear functions defined by the map $\phi_{r'} \colon N_{r'}
\rightarrow L_r \isom \RR^{d+n-1}$ from Definition~\ref{d:linear}. First suppose
that $r'$ intersects the interior of a facet~$f$ of~$\Delta$. Then, $N_{r'}$ is
the union of the two facets containing $r'$, and naturally embeds in $\RR^{2+n}$
by sending each vertex to a coordinate vector, as in Definition~\ref{d:linear}.
The inclusion of $N_{r'}$ into $f \subset \RR^n$ factors through a linear map
from $\RR^{2+n}$ to $\RR^n$, whose kernel contains $(1, 1, -\alpha(v_1', r'),
\ldots, -\alpha(v_n', r'))$ by the choice of those structure constants in
Construction~\ref{const:subdivision}. Therefore, the inclusion of $N_{r'}$ into
$\RR^n$ factors through the map $\phi_{r'} \colon N_{r'} \rightarrow L_{r'}$,
defining linear functions on the open subset $N_{r'}^o$, and so the linear
functions on $\Delta$ and $\Delta'$ agree.

Now suppose that $r'$ is contained in a ridge $r$ of $\Delta$, and
we let $v_i$, $v_i'$, $x_{i,j}$ and $y_{i,j}$ be as
in Construction~\ref{const:subdivision}, and $\phi_r$ as in
Definition~\ref{d:linear}. Then, using the fact that $\phi_r$ is affine linear
on each facet of $\Delta$, we can evaluate:
\begin{align*}
\sum_{i=1}^d \phi_r(w_i') &=
\sum_{i=1}^d \frac{1}{m} \Bigg(\phi_r(w_i) + \sum_{j = 0}^n x_{i,j}
\phi_r(v_j)\Bigg) \\
&= \frac{1}{m}\sum_{j=0}^n \alpha(v_j, r) \phi_r(v_j) +
\frac{1}{m}\sum_{i=1}^d\sum_{j=0}^n x_{i,j}\phi_r(v_j) \\
&= \frac{1}{m} \sum_{j=0}^n \Bigg(\alpha(v_j, r) + \sum_{i=1}^d x_{i,j} \Bigg)
\phi_r(v_j) \\
&= \sum_{j=0}^n \alpha(v_j', r') \phi_r(v_{j}'),
\end{align*}
where the last step is by the change of variables in~(\ref{eq:subdivision}),
together with the linearity of $\phi_r$ on~$r$. Thus, $\phi_r$ satisfies
the defining linear relation of~$\phi_{r'}$, and so $\phi_{r'}$ is isomorphic to
the restriction of $\phi_r$, as in the previous paragraph.
\end{proof}

If $\Delta'$ is a subdivision of a weak tropical complex $\Delta$, then any
formal sum of polyhedra on $\Delta'$ is also a formal sum of polyhedra on
$\Delta$. Conversely, for any formal sum of polyhedra on $\Delta$, we can
intersect with all the facets of $\Delta'$ to obtain an equivalent formal sum
of polyhedra, all of whose terms are contained in a single facet of $\Delta'$,
and thus can be considered as a formal sum of polyhedra on $\Delta'$. Moreover,
this correspondence passes to equivalence class of formal sums of polyhedra,
so we can identify the equivalence classes on $\Delta$ with those on $\Delta'$,
and by the following lemma, we can also identify divisors on $\Delta$ with those
on~$\Delta'$.

\begin{lem}\label{l:subdivide-divisors}
Let $\Delta'$ be an order~$m$ subdivision of a weak tropical complex~$\Delta$.
Then a formal sum of polyhedra is a divisor on $\Delta$ if and only if it is a
divisor on~$\Delta'$.
\end{lem}

\begin{proof}
By Lemma~\ref{l:equiv-balanced}, passing to equivalent formal sums of polyhedra
doesn't change whether or not it is balanced, so we can assume that the terms of
the formal sum of polyhedra are each contained in the facets of~$\Delta'$. Then,
Lemma~\ref{l:subdivide-linear} shows the linear functions on $\Delta$ and
$\Delta'$ are the same, up to multiplication by $m$, so we just need to check
that the multiplicity $\mult_{Z,Q}(\phi)$ of a linear function along $Q$ of a
formal sum of polyhedra $Z$ vanishes in $\Delta$ if and only if it vanishes
in~$\Delta'$.

Here, the key the difference is between two different embeddings of a single
facet $f'$ of $\Delta'$ into $\RR^n$. First, we can identify $f'$ with a
standard unit simplex in $\RR^n$. Second, we identify $f'$ with a simplex
contained in a facet $f$ of $\Delta$, and $f$ is identified with a standard unit
simplex in $\RR^n$. Because the subdivision of $\Delta$ is assumed to be
unimodular in Construction~\ref{const:subdivision}, in the latter
identification, $f'$ is a standard unit simplex scaled by $1/m$. Thus, if we
transform the vectors $\mathbf w_i$ and $\mathbf v_i$ used to compute the
multiplicity in $\Delta$ into the coordinates where $f$ is identified with a
standard unit simplex, then the transformed vectors have a common multiple of
$m$ in their coordinates. Scaling these vectors by $1/m$ gives suitable vectors
for computing the multiplicity in $\Delta'$, and so the multiplicity of a given
linear function along $Q$ of $Z$ in $\Delta'$ is $1/m$ times the multiplicity in
$\Delta$, which completes the proof.
\end{proof}

\begin{lem}\label{l:subdivide-lin-equiv}
Let $\Delta'$ be an order~$m$ subdivision of a weak tropical complex~$\Delta$,
and further assume that every ridge of $\Delta$ is contained in a facet.
Then, divisors on $\Delta$ are linearly equivalent if
and only if they are linearly equivalent on~$\Delta'$
\end{lem}

\begin{proof}
To show that linear equivalence on $\Delta$ agrees with that on $\Delta'$,
we need to show that, if $\phi$ is a PL function on $\Delta$, then
$\div(\phi) = \div'(m \phi)$, where $\div$ and $\div'$ are the divisor functions
from Proposition~\ref{p:div}, on $\Delta$ and~$\Delta'$, respectively. Since
Proposition~\ref{p:div} uniquely characterize $\div$ and $\div'$ in terms of
four properties, we need to check that those properties coincide for $\Delta$
and~$\Delta'$. The first two properties are that the divisor is linear in the PL
function and local, which are independent of the subdivision. Moreover, the
third property is that linear functions are characterized by having trivial
divisor, and by Lemma~\ref{l:subdivide-linear}, the linear functions on~$\Delta$
and~$\Delta'$ agree, up to scaling. The fourth property normalizes the divisor
function via a PL function supported on a single facet, which also agrees,
because when passing to~$\Delta'$ we scale both the coordinates on the complex
and the PL function by the same factor of~$m$.
\end{proof}

On the other hand, the subdivisions of Construction~\ref{const:subdivision}
mirror what happens for toroidal resolutions of ramified base changes.

\begin{lem}\label{l:subdivide-toroidal}
Let $\X$ be a regular semistable degeneration with dual complex~$\Delta$. If
$\Delta'$ is an order~$m$ subdivision of~$\Delta$ and $R'$ is a totally ramified
degree~$m$ extension of~$R$, then the corresponding resolution~$\X'$ of
$\X \times_{\Spec R} \Spec R'$ has $\Delta'$ as its weak tropical complex.
\end{lem}

\begin{proof}
We first describe the pullback of a divisor~$C_v$ from $\X$ to~$\X'$. If $v'$ is
a vertex of~$\Delta'$, contained in a $k$-dimensional simplex~$s$ of~$\Delta$,
then we can express the coordinate of $v'$ by a vector $(x_0, \ldots, x_k)$ with
$x_0 + \cdots + x_k = m$. If the vertices of $s$ are $v_0, \ldots, v_k$, then we
claim that the coefficient of $C_{v'}$ in the pullback of $C_{v_i}$ is $x_i$.
The reason is that in local, toroidal coordinates around $C_{s}$, the divisor
$C_{v_i}$ is defined by a monomial with exponent vector $(0, \ldots, 1, \ldots,
0)$. The valuation of the divisor $C_v'$ in $\X'$ corresponds to the vector
$(x_0, \ldots, x_k)$ in the lattice dual to these monomials, meaning that the
valuation of the monomial locally defining $C_{v_i}$ is given by the dot
product of these two vectors, which is $x_i$.

We now adopt the notation of Construction~\ref{const:subdivision}. For any facet
$f'$ containing $v'$, we compute
the intersection number
$\pi^{-1}(C_{v_i}) \cdot C_{f'}$ in two different ways. First, using the
discussion in the previous paragraph, we can represent $\pi^{-1}(C_{v_i})$ as a
linear combination of divisors in~$\X'$. The ones which intersect $C_{f'}$
correspond to the vertices of~$f'$, together with the vertices in neighboring
simplices. All together, these give:
\begin{equation*}
x_{1,i} + \cdots + x_{d,i} -
\alpha(v_0', f') y_{0,i} - \cdots - \alpha(v_n', f') y_{n,i}
\end{equation*}
On the other hand, using the projection formula, $\pi^{-1}(C_{v_i}) \cdot
C_{f'} = C_{v_i} \cdot \pi(C_{f'})$, and we know that $\pi_*(C_{f'})$ is $C_f$,
so the intersection number is $-\alpha(v_i, f)$. Putting these equations
together for all $i$, we get:
\begin{equation*}
\begin{pmatrix}
y_{0,0} & \cdots & y_{n,0} \\
\vdots & & \vdots \\
y_{0,n} & \cdots & y_{n,n}
\end{pmatrix}
\begin{pmatrix} -\alpha(v_0', F') \\ \vdots \\ -\alpha(v_n', F') \end{pmatrix}
+ \begin{pmatrix}
x_{1,0} + \cdots + x_{d,0} \\
\vdots \\
x_{1,n} + \cdots + x_{d,n}
\end{pmatrix}
=
\begin{pmatrix} -\alpha(v_0, F) \\ \vdots \\ -\alpha(v_n, F) \end{pmatrix}.
\end{equation*}
We can solve for the $\alpha(v_i', F')$ and we get~(\ref{eq:subdivision}).
\end{proof}

\section{Refined specialization}\label{s:refined}

In this section, we define a refined specialization map from divisors on the
general fiber of~$\X$ to divisors on the weak tropical complex~$\Delta$. Given a
divisor~$D$, we first choose a ramified extension of the discrete valuation
ring~$R$, followed by a toroidal resolution of singularities~$\X'$ such that if
$D'$ is the pullback of $D$ to the general fiber of $\X'$, the closure
$\overline D'$ does not contain any strata of the special fiber of~$\X'$. Then,
the refined specialization of~$D$ is $\rho(D')$, which can be considered as a
divisor either on the subdivision~$\Delta'$ or on~$\Delta$ by the results in the
previous section. The refined specialization is the higher-dimensional analogue
of the map~$\tau_*$ from~\cite{baker}*{Sec.~2C}. As stated in the introduction,
the purpose of the refined specialization is that the refined specialization of
an effective divisor is always effective, which is not necessarily true for the
coarse specialization~$\rho$.

\begin{figure}
\includegraphics{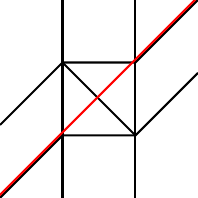}
\caption{The black lines give a subdivision of the plane giving a semistable
family whose general fiber is $\PP^1 \times \PP^1$. The red line is the
tropicalization of a line in one of the rulings whose specialization to the
tropical complex consists the sum of the outside edges of the square minus the
diagonal edge. In particular, the specialization is not effective.}
\label{f:non-effective}
\end{figure}

\begin{ex}\label{ex:non-effective}
We construct our example as a toric degeneration. In particular, we take the
3-dimensional fan whose support is $\RR^2 \times \RR_{\geq 0}$ and intersection
with the plane $\RR^2 \times \{1\}$ is the configuration in
Figure~\ref{f:non-effective}. Projection to the last factor defines a toric
morphism from a 3-dimensional toric variety $\X$ to $\mathbb A^1$, whose fibers
are isomorphic to $\PP^1 \times \PP^1$, except over the origin. The fiber of the
origin consists of 4 components, corresponding to the 4 rays generated by the 4
vertices in the Figure~\ref{f:non-effective}. The intersections between these
components correspond to the cones they have in common, such that the dual
complex of $\X$ is the triangulation of a square, making up the bounded cells in
Figure~\ref{f:non-effective}.

By Proposition~\ref{p:robust-dim-1}, the locally closed strata of dimension 1 in
$\X$ are affine and we claim that the same is true in dimension 2. Each of the
four components is isomorphic to the toric variety corresponding to the star of
the corresponding vertex~$v$, with the divisor $D_v$ equal to the union of the
rays corresponding to the bounded edges. Therefore, in each case, the locally
closed stratum $C_v \setminus D_v$ is isomorphic to $\mathbb A^2$.

Now consider a line in one of the rulings of $\PP^1 \times \PP^1$ whose
tropicalization is the red line in Figure~\ref{f:non-effective}. The
specialization of this line to~$\Delta$ is the sum of four outside edges of the
square minus the inner diagonal.  In particular, this specialization is not
effective and it does not agree with the tropicalization. However, the
specialization is linearly equivalent to an effective divisor, namely the
intersection of the red tropicalization with~$\Delta$. 
\end{ex}

In Example~\ref{ex:non-effective}, it was easy to guess the refined
tropicalization. In general, we need to show that there exists a subdivision of
the weak tropical complex~$\Delta$ which will give the refined specialization
according to the definition at the beginning of this section, which is the
content of the following lemma. The construction is essentially a toroidal
version of the partial resolution of a hypersurface singularity coming from its
Newton polyhedron.

\begin{lem}\label{l:fine-effective}
Let $D$ be an effective divisor on the general fiber~$X$ of a
degeneration~$\X$. Then, there exists a finite ramified extension~$R'$ of~$R$
and a toroidal resolution of singularities $\X'$ of $\X \otimes_R R'$ such that
$\overline D'$ doesn't contain any of the strata of $\X'$, where $\overline D'$
is the closure in~$\X'$ of the pullback of~$D$ to the general fiber~$X'$.
\end{lem}

\begin{proof}
We let $s$ be any $k$-dimensional simplex of~$\Delta$. Let $A$ be the local
ring of $\X$ at the generic point of~$C_s$. By assumption, $A$ is a regular
ring with maximal ideal generated by $x_0, \ldots, x_k$, each of which is the
local
equation of one of the components of~$\X_0$ containing~$C_s$. The
completion~$\widehat A$ of~$A$ is isomorphic to $S[[x_0, \ldots, x_k]]/\langle
x_0 \cdots x_k - \pi\rangle$, where $S$ is a complete, local extension of~$R$
and $\pi$ is a uniformizer in~$S$. We let $f \in \widehat A$ be the local
equation of $\overline D$ in the completed local ring~$\widehat A$. We construct
the local Newton polyhedron from~$f$ as follows. For each term $s_a x_0^{a_0}
\cdots x_k^{a_k}$ in the summation representation of~$f$, we translate the
positive orthant in~$\RR^{k+1}$ by $(a_0 + \val(s_a), \ldots, a_k + \val(s_a))$.
Then, the local Newton polyhedron of~$f$ is the convex hull of the translated
orthants coming from all terms of~$f$. If $f$ is not a unit, this convex hull
will not contain the origin. In this case, we take the dual polyhedron of the
local Newton polyhedron and project its bounded faces to the 
unit $k$-dimensional simplex along lines
through the origin. This yields a rational polyhedral subdivision of~$s$.

We now claim that these subdivisions are compatible so that taken together they
form a rational polyhedral subdivision of~$\Delta$. It is sufficient to show
compatibility between the subdivision of~$s$ and of one of its faces~$s'$, which
we can also take to be the face indexed by~$k$.
Thus, in the local ring~$A$, the component $C_{s'}$ is
defined by the ideal $\langle x_0, \ldots, x_{k-1}\rangle$, and so its local
ring is the
localization of~$A$ at this ideal. We consider the image of~$f$, the local
equation of~$D$ in the tensor product $\widehat A \otimes_A A_{\langle x_0,
\ldots, x_{k-1}\rangle}$. We can now regroup the terms of $f$ as:
\begin{equation*}
\sum_{a_0, \ldots, a_{k-1}} \bigg(\sum_{a_k} s_ax_k^{a_k}\bigg) x_0^{a_0} \cdots
x_{k-1}^{a_{k-1}}.
\end{equation*}
This expression also holds in the completion of the local ring of $C_{s'}$,
so the local Newton polyhedron of $f$ along $C_{s'}$ is
the projection of the polyhedron along $C_s$. Thus the dual of the Newton
polyhedron is the restriction of the dual, which shows the desired
compatibility.

Now we have a polyhedral subdivision of~$\Delta$.
If $m$ is the least common multiple of the denominators of the vertices of the
subdivision, then it corresponds to a ramified degree~$m$ extension $R'' \supset
R$ followed by a toroidal map $\X'' \rightarrow \X \times_R R''$. Although
$\X''$ is not semistable, the closure $\overline D''$ of the pullback of $D$ to
$\X''$ does not contain any of the strata. By~\cite{kkmsd}*{Thm.~III.4.1}, there
exists a further extension and toroidal map $\X' \rightarrow \X'' \times_{R''}
R'$ such that $\X'$ is semistable. Since each stratum of $\X'$ maps to a stratum
of $\X''$, the pullback of $\overline D''$ to~$\X'$ doesn't contain any strata
of the special fiber either. This implies that the pullback of $\overline D''$
is equal to the closure of the pullback of $D$.
\end{proof}

After the subdivision given by Lemma~\ref{l:fine-effective}, the
specialization of~$D$ will be effective, because all of the intersections
between $\overline D$ and curves $C_r$ will be proper, and thus the coefficients
making up the definition of the specialization will all be non-negative. In
addition to the above  construction, the refined specialization has a more
intrinsic interpretation, at least set-theoretically, as the projection of the
Berkovich analytification of the divisor to the skeleton defined by the
degeneration. Here, we are using the deformation retract from the Berkovich
analytification~$X^{\operatorname{an}}$ to the dual complex of~$\X$. Such
a retract is defined in a more general setting in~\cite{berkovich2} and in the
specific case of regular semistable degenerations in~\cite{nicaise}*{Sec.~3.3}.

\begin{proof}[Proof of Proposition \ref{p:analytification}]
By Lemma~\ref{l:fine-effective}, we can assume that our degeneration~$\X$ is
chosen such that $\overline D$ doesn't contain any of the closed strata of~$\X$.
Because projection of $X^{\operatorname{an}}$
to the skeleton is compatible with the specialization map, we know that the
image of~$D^{\operatorname{an}}$ is contained in those simplices~$s$ such that $\overline D$
intersects~$C_s$. By hypothesis, $\overline D$ doesn't intersect any
$0$-dimensional strata
$C_f$ and whenever it intersects a curve~$C_r$ for some ridge~$r$, it does so
properly and thus $r$ has positive coefficient in the specialization of~$D$.
Therefore, the projection of the analytification is contained in the union of
the support of the refines specialization of~$D$,
together with some simplices of dimension $0 \leq n-2$.

It now only remains to show that the projection map is surjective onto the
ridges~$r$ in the support of the specialization. For this we suppose that $r$ is
a ridge such that $\overline D$ intersects~$C_r$ and let $x$ be a point in the
intersection. In a neighborhood of~$x$, we can take local defining equations for
the divisors containing~$C_r$ and use them to define a map to $\Spec S$ where $S
= R[x_1, \ldots, x_{n}] / \langle x_1 \cdots x_{n} - \pi \rangle$. As
in~\cite{nicaise}*{Sec.~3.3, Case~1}, the $(n-1)$-dimensional simplex, which we
identify with $r$ embeds in the analytification of $\Spec S$ as the skeleton.
Since the fiber over $x_1 = \cdots = x_n = 0$ is finite, the map from $\overline
D$ to $\Spec S$ is dominant and so by Proposition~3.4.6(7) of~\cite{berkovich},
we can lift any norm corresponding to a point in this skeleton to a point in the
analytification of $D$. Since the projection onto the skeleton of~$X$ is
defined in terms of the norms of the~$x_i$, we've produced a norm in the
analytification of~$\overline D$ which projects onto any point of~$r$.
\end{proof}

The image of the projection in Proposition~\ref{p:analytification} can be
larger than the refined specialization whenever $\overline D$ meets a
component~$C_v$, but doesn't meet any curves~$C_r$ contained in~$C_v$.

\begin{ex}\label{ex:blow-up}
Suppose $\X$ is any degeneration of dimension at least~$2$ and $x$ is a
point of $X$ whose closure in~$\X$ is contained in $C_v \setminus D_v$ for
some vertex~$v$. Then, the blow-up~$\X'$ of~$\X$ at the closure of~$x$ is a
regular semistable degeneration and the specialization of the exceptional
divisor~$E$ of the blow-up is trivial because its closure doesn't intersect any
curve $C_r$ in~$X$. However, the analytification of $E$ is not
empty and by the compatibility with specialization noted in the proof of
Proposition~\ref{p:analytification}, its projection to~$\Delta$ must be the
single point~$v$, which has codimension $n \geq 2$.
\end{ex}

\section{Proof of the specialization inequality}\label{s:proof}

In this section, we finish the proof of Theorem~\ref{t:specialization}. The
main remaining ingredient is to show that the refined specialization always
preserves point containment, at least for appropriately chosen points.
Example~\ref{ex:blow-up} shows that even for robust degenerations, the
specialization of a non-trivial divisor can be trivial, which the
importance of the choice of a point. However, we have the following
result with the locally closed strata are affine.

\begin{prop}\label{p:locally-closed-affine}
Suppose that $\X$ is a degeneration whose locally closed strata of dimension~$m$
are affine for all
$2 \leq k \leq m$. If $E$ is a divisor in~$X$ and its closure~$\overline
E$ intersects some stratum~$C_s$ of the special fiber of~$\X$ and $C_s$ has
dimension~$m$, then $\overline E$ intersects the curve~$C_r$ for some ridge~$r$
containing~$s$.
\end{prop}

\begin{proof}
The proof is by induction the dimension~$m$. If $m = 1$, then $s$ is a ridge
and we're done. Otherwise, by assumption, $C_s \setminus D_s$ is affine,
and affine varieties contain no complete subvarieties of positive dimension.
However, $\overline E \cap C_s$ is a divisor in $C_s$ and so of dimension $m -
1 \geq 1$. Thus, $\overline E$ must intersect $D_s$ and thus $C_{s'}$ for some
$(n-m+1)$-simplex $s'$ containing $s$. Applying the inductive hypothesis, we
get that $\overline E$ intersects $C_r$ for some ridge~$r$.
\end{proof}

Proposition~\ref{p:locally-closed-affine} gives us a tool for proving that
the refined specialization contains a given vertex of the dual complex, and
together with Lemma~\ref{l:subdivision-refined} below, it can be used to give
a proof of Theorem~\ref{t:specialization} for degenerations whose locally
closed strata are affine. However, Theorem~\ref{t:specialization} applies not
only to degenerations where the locally closed strata are affine, but also to
those which are only robust in the top two dimensions. In such cases, the
conclusion of Proposition~\ref{p:locally-closed-affine} doesn't necessarily
hold, as shown by Example~\ref{ex:blow-up}. However, we can still ensure that
the refined specialization preserves point containment, so long as the point is
chosen sufficiently generically.

\begin{lem}\label{l:lifting-point}
Assume that $R$ is complete, with fraction field $K$, and that $\X$ is robust in
dimension $n-1$ and~$n$ and that the locally closed strata of dimension at most
$n-2$ are affine. Then for any vertex~$v$ in~$\Delta$, there exists a
$K$-point~$x$ in~$\X$ such that for any effective divisor~$E \subset X$
containing~$x$, the refined specialization of~$E$ contains~$v$.
\end{lem}

Since the crux of this lemma is using the robustness hypothesis, we first prove
a couple of results about big divisors. We begin with an application of
Bertini's theorem to the connectedness of big divisors.

\begin{lem}\label{l:big-connected}
Any big Cartier divisor has a unique connected component which is itself a big
Cartier divisor.
\end{lem}

\begin{proof}
Suppose that $E$ is a big divisor on a variety~$X$. Then we let $\widetilde X$
be the blow-up of the base locus and let $\mathcal O_{\widetilde X}(1)$ be the
relative ample divisor. Let $\widetilde E$ be the pullback of~$E$ and then
$\widetilde E(-1)$ defines a regular morphism from~$\widetilde X$ to~$\PP^N$.
By Bertini's theorem~\cite{jouanolou}*{Thm.~7.1}, the support of~$\widetilde E(-1)$ is connected.
Therefore, the base locus of~$E$ contains all but one connected component
of~$E$. By taking just that connected component, we've only removed fixed
divisors, and so we still have a big divisor, which is the desired statement.
\end{proof}

\begin{lem}\label{l:contract-big}
Suppose that $E$ is a big Cartier divisor on a variety~$X$. If $\phi$ is any
morphism from $X$ to~$\PP^N$, then either $\phi(E) = \phi(X)$ or $\phi$ is
generically finite on some irreducible component of~$E$.
\end{lem}

\begin{proof}
Suppose, for contradiction, that $E$ is a big divisor on~$X$ and $\phi\colon X
\rightarrow \PP^N$ is a morphism which is not generically finite on any
component of~$E$ and such that $\phi(X) \neq \phi(E)$. We first
consider the case that $\phi$ is generically finite on~$X$. Then, $\dim(X) =
\dim\phi(X) \geq \dim \phi(E) + 2$, so for any point in $\phi(X) \setminus
\phi(E)$, we can take an intersection with a general linear space to obtain
a curve through $\phi(X)$ which doesn't intersect
$\phi(E)$. By taking the preimage, we get a curve passing through any
sufficiently general point in~$X$, which does not intersect~$E$. Since $E$ is
big, such curves can't exist, and so we get a contradiction.

Second, we suppose that $\phi$ is not generically finite on~$X$, but still that
$\phi(X) \neq \phi(E)$. Then for any point in $\phi(X) \setminus \phi(E)$, the
fiber of~$\phi$ is at least $1$-dimensional and we can choose a curve in this
fiber. Again, we've constructed curves in~$X$ which pass through sufficiently
general points, and which do not intersect~$E$, so we've again contradicted our
assumption that $E$ is big.
\end{proof}

\begin{proof}[Proof of Lemma~\ref{l:lifting-point}]
Because $R/\m$ is algebraically closed, we can choose a general point~$x_0$
in~$C_v \setminus D_v$. By general, we mean that it lies outside of finitely
many closed subvarieties, which don't depend on~$E$ and will be made explicit in
the course of the proof. By Hensel's lemma, we can lift the point~$x_0$ to an
$R$-point of~$\X$, which will give the desired point~$x$.

Now let $E$ be any divisor of $X$ containing $x$ and we first want to show
that its closure $\overline E$ in $\X$ intersects a $1$-dimensional
stratum~$C_r$ for some ridge~$r$ containing $v$. If the dimension~$n$ is~$1$, then this is immediate, so we assume
that $n \geq 2$. By assumption, some multiple of~$D_v$ defines a rational map
$\phi_v \colon C_v \dashrightarrow \PP^N$ which is birational onto its image. We
can assume that our choice of~$x_0$ was in the locus where $\phi_v$ defines an
isomorphism onto its image. Since $\overline E$ contains $x_0$, the restriction
of~$\phi_v$ to~$E_0 = \overline E \cap C_v$ must also be birational onto its
image. Since this restriction is defined by the pullback of the same multiple of
$D_v$, we see that $D_v \cap E_0$ is a big divisor on~$E_0$.

By Lemma~\ref{l:big-connected}, we can take $B$ to be a connected component of
$D_v \cap E_0$ which is a big divisor on~$E_0$. Since $B$ is connected, it is
either contained in a single stratum~$C_e$ for some edge~$e$ or it meets some $C_t$,
where $t$ is a $2$-simplex of~$\Delta$. If $n=2$, then $C_e$ is a
$1$-dimensional stratum intersecting $\overline E$, which is what we wanted to
show. So, we assume that $n \geq 3$ and, for the moment, we also assume that $B$
is contained in a single~$C_e$.

By assumption, $D_e$ is a big divisor in $C_e$, so some multiple of~$D_e$
defines a map $\phi_e \colon C_e \dashrightarrow \PP^M$, birational onto its
image. Since this
morphism is just defined by a collection of rational functions on $C_e$, we can
extend it to a rational map from~$C_v$ to~$\PP^M$, also denoted
by~$\phi_e$. The image~$\phi_e(C_v)$ contains the image of~$C_e$, so it is at
least $(n-1)$-dimensional. Therefore, a general fiber of $\phi_e$ is at most
$1$-dimensional, so we can assume that $\phi_e^{-1}(\phi_e(x_0))$, the fiber
containing~$x_0$, has dimension at most~$1$. Thus, $\phi_e(E_0)$ must be at
least $(n-2)$-dimensional. We apply Lemma~\ref{l:contract-big} to the
restriction~$\phi_e\vert_{E_0}$, and in either of that lemma's two
cases, $\phi_e(B)$ is also at least $(n-2)$-dimensional. Since $\phi_e(B)$ is
positive dimensional, it intersects any hyperplane in $\PP^M$, which shows that
$B$ must intersect $D_e$, and thus $C_t$ for some $2$-simplex~$t$.

At this point, we've shown that if $n \leq 3$, then $\overline E$ intersects
$C_r$ for some ridge~$r$ containing~$v$. If $n > 3$, then $\overline E$ at least
intersects $C_t$ for some $2$-simplex~$t$ containing~$v$, but we can then apply
Proposition~\ref{p:locally-closed-affine} to show that $\overline E$
intersects some $C_r$ for some ridge containing~$v$ in this case as well.

We now consider the resolution of the base change~$\X'$ produced by
Lemma~\ref{l:fine-effective}. There are two cases, depending on whether or not
there exists a facet~$f$ containing~$r$ such that $\overline E$ contains the
point~$C_f$. If there does, then we consider the component $C_v'$ in~$\X'$
corresponding to~$v$, which maps birationally onto $C_v$. By the properness of
the resolution, there must be a point in $C_v' \cap \overline E'$ mapping to the
point $C_f$, where $\overline E'$ is the closure of $E$ in~$\X'$. However, the
fiber of $C_f$ in the resolution~$\X'$ is a union of toric varieties
intersecting along their toric boundaries and so by
Proposition~\ref{p:locally-closed-affine}, $\overline E'$ must intersect one
of the $1$-dimensional strata in this fiber. By Lemma~\ref{l:fine-effective},
this intersection is proper, so it gives a positive coefficient to some ridge
containing~$v$.

On the other hand, suppose that $\overline E$ does not intersect $C_f$ for
any facet~$f$ containing~$r$. We've assumed that $D_r$ is a big divisor
on~$C_r$, so it is a non-trivial divisor. In other words, there is at least one
stratum $C_f$ properly contained in $C_r$. Since $\overline E$ does not
contain any $C_f$, it must intersect $C_r$ properly, and thus
with positive intersection number. Now we choose a ridge~$r'$ containing $v$ in
the subdivision of~$r$ and then the corresponding stratum $C_{r'}'$ of~$\X'$
maps surjectively onto $C_r$. By the projection formula, this implies that the
intersection number of $\pi^{-1}(\overline E)$ with $C_{r'}'$ is positive, where
$\pi$ is the map from~$\X'$ to~$\X$. Of course, $\pi^{-1}(\overline E)$ is not
the same as~$\overline E'$, but also includes any
components of the special fiber whose image is contained in $\overline E$.
However, the only such components which affect the intersection number will be
those which meet $C_{r'}'$, in which case the projection of the component onto
$\X$ must be either $C_r$ or a point $C_f$ contained in $C_r$. However, we've
assumed that in neither of these can be contained in~$\overline E$, so
$\overline E' \cdot C_{r'}' = \pi^{-1}(\overline E) \cdot C_{r'}'$ is positive,
which shows that the ridge $r'$ containing~$v$ occurs with positive coefficient
in the refined specialization, as desired.
\end{proof}

We now have one lemma remaining before finishing the proof of
Theorem~\ref{t:specialization}. While
Lemma~\ref{l:lifting-point} shows how to find a divisor whose refined
specialization contains a given vertex, our definition of~$h^0$ allowed
arbitrary rational points. We can turn arbitrary rational points into vertices
by choosing an appropriate subdivision, but we need to check that such a
subdivision preserves the robustness and affine properties of the
degeneration.

\begin{lem}\label{l:subdivision-refined}
Let $\X$ be a tropical complex and let $\X'$ be any toroidal resolution of a
base change of~$\X$ to a totally ramified extension of~$R$.
Let $m$ be an integer. If the $k$-dimensional
locally closed strata of~$\X$ are affine for all $k \leq m$,
then the same is true
for the $k$-dimensional locally closed strata in $\X'$ when $k \leq m$.
If the $\X$ is robust in dimension~$k$ for all $k \leq m$,
then $\X'$ is also robust in dimensions $k \leq m$.
\end{lem}

\begin{proof}
Let $\X'$ be a toroidal resolution of a ramified extension of~$\X$.
We first consider the case that the $k$-dimensional locally closed strata of
$\X$ are affine for
$k \leq m$. Suppose that $C_{s'}$ is a stratum of~$\X'$ which has
dimension $k' \leq m$. Then the image of $C_{s'}$ in $\X$ is the stratum $C_s$,
where $s$ is the minimal simplex of $\Delta$ which contains $s'$. Therefore,
$k$, the dimension of $C_s$ satisfies $k \leq k' \leq m$, and so $C_s \setminus
D_s$ is affine by assumption. From the construction of the toroidal resolution,
we know that $C_{s'}$ will be a toric variety bundle over $C_s$ and $C_{s'}
\setminus D_{s'}$ will be a $\Gm^{k' - k}$-bundle over $C_s \setminus D_s$.
Therefore, the map from $C_{s'} \setminus D_{s'}$ to $C_s \setminus D_s$ is
affine and so $C_{s'} \setminus D_{s'}$ is affine as well. We conclude that
the $k$-dimensional locally closed strata of~$\X'$ are affine for $k \leq m$.

We now consider the case that $\X$ is robust in dimensions $k$ for $k \leq m$.
Let $C_{s'}$ and $C_s$ be strata of $\X'$ and~$\X$ of dimensions $k'$ and $k$ as
in the previous paragraph. By assumption, some multiple of the divisor~$D_s$
defines a birational map $\phi$ from an open subset of~$C_s$ to $\PP^N$. We let
$\overline{\phi(C_s)}$ denote the closure of the image of this map and let $Y$
denote the image of $C_s \setminus D_s$. Since $Y$ is the complement of a
hyperplane section in a projective variety, $Y$ is affine.

Now, we consider $C_{s'}$, which, as above, is a toric variety bundle
over~$C_s$, and the toric variety is described by the star of~$s'$ in~$s$. We
can choose some positive multiple of the boundary of this toric variety which
contains a spanning set of the characters of the torus, i.e.\ defines a
birational map to projective space. Let $\mathcal L$ be the corresponding line
bundle on~$C_{s'}$. Explicitly, $\mathcal L$ is the line bundle associated to a
multiple of the sum of the divisors corresponding to the simplices of~$\Delta'$
containing~$s'$, which are also contained in~$s$. We consider the coherent sheaf
on~$Y$, $(\phi\circ\pi)_* (\mathcal L)$, where $\pi \colon C_{s'} \rightarrow
C_s$ is the restriction of the map from $\X'$ to~$\X$.

We can find an open set~$U$ of~$Y$ such that $\phi$ is an isomorphism on~$Y$ and
$(\phi \circ \pi)^{-1}(U)$ is the trivial toric variety bundle over
$\phi^{-1}(U) \isom U$. Furthermore, we can assume that $U$ is the complement of
the variety defined by some element~$f$ of the global sections of~$Y$, since $Y$
is affine. On~$U$, the push-forward $(\phi \circ \pi)_*(\mathcal L)$ is a free
sheaf whose generators can be identified with the sections on the toric variety.
In particular, these sections define a rational map from $(\phi \circ
\pi)^{-1}(U)$ to projective space which is birational on each fiber of~$\pi$.
Since $Y$ is affine, we lift these sections from~$U$ to~$Y$ after multiplying by
a sufficiently large power of~$f$. If we regard $f$ as a section of $\mathcal
O(\ell D_s)$ for sufficiently large $\ell$, then what we've found are sections of
$\mathcal L \otimes \pi^{-1}(\mathcal O(\ell D_s))$ which define an embedding of
the torus for a generic fiber of~$\pi$. Combining these with the pullbacks of
the sections defining~$\phi$, we get a birational map from $C_{s'}$, and
therefore $D_{s'}$ is a big divisor and $\X'$ is robust in dimensions~$k \leq
m$.
\end{proof}

\begin{proof}[Proof of Theorem~\ref{t:specialization}]
By base changing, we can assume that $R$ is complete without changing
$h^0(X, D)$ or~$\Delta$. Let $r$ be an integer less than $\dim
H^0(X, \mathcal O(D))$ and let $p_1, \ldots, p_r$ be any $r$ rational
points in~$\Delta$. We want to show that there exists an effective divisor
linearly equivalent to $\rho(D)$ containing these points.

We first make a subdivision of $\Delta$ as follows. Let $m$ be the least common
denominator of all the coordinates of the points~$p_i$. We rescale the simplices
by~$m$ and subdivide them using the ``regular subdivision''
of~\cite{kkmsd}*{Thm.~III.2.22}, and now the points~$p_i$ are vertices. We let
$\X'$ be the corresponding toroidal resolution of a ramified extension of~$R$.
By Lemma~\ref{l:subdivision-refined}, the locally closed strata of $\X'$ of
dimension at most $n-2$ are affine, and it is robust in dimensions $n$ and
$n-1$. By Lemma~\ref{l:subdivide-toroidal}, the tropical complex of $\X'$ is the
subdivision of $\Delta$ as in Construction~\ref{const:subdivision}. Finally, by
Lemmas~\ref{l:subdivide-divisors} and~\ref{l:subdivide-lin-equiv}, if we can
find a divisor linearly equivalent to $\rho(D)$ on~$\Delta'$, then it will also
be a divisor linearly equivalent to $\rho(D)$ on $\Delta$. Therefore, we replace
$\X$ and $\Delta$ with $\X'$ and $\Delta'$ for the rest of the proof.

Using Lemma~\ref{l:lifting-point}, we can choose points $x_1, \ldots, x_r$
in~$X$ corresponding to $p_1, \ldots, p_r$ respectively. Since each $x_i$ is
rational over the fraction field of $R$, vanishing on~$x_i$ imposes one linear
condition on the sections $H^0(X, \mathcal O(D))$. We've assumed that $r$ is
less than the dimension of this vector space, so there exists a non-zero section
of $H^0(X, \mathcal O(D))$ defining a divisor $D'$ which contains all of
the~$x_i$. Let $D''$ be the refined specialization of~$D'$. Then $D''$ is
effective by Lemma~\ref{l:fine-effective} and contains the points $p_1, \ldots,
p_r$ by Lemma~\ref{l:lifting-point}. By \cite{cartwright-complexes}*{Prop. 5.2},
$D''$ is linearly equivalent to $\rho(D)$. Therefore, $h^0(\Delta, \rho(D))$ is
greater than $r$, which establishes the desired inequality.
\end{proof}

\end{document}